\documentclass[11pt]{article}
\pdfoutput=1
\usepackage{amsmath,amssymb,amsthm}
\usepackage{gbmacros}
\usepackage{times,mathpazo,charter,graphicx,subfigure,ulem}
\usepackage{hyperref,cite}

\makeatletter
\newcommand*\titlesize{\@setfontsize\titlesize{16}{1.2}}
\makeatother
\def\title#1{{\noindent\titlesize\bf\sffamily#1\par}\vglue2ex}
\def\author#1{\par{\noindent\sffamily\large#1}\par\vglue1.4ex}
\def\address#1{{\it\noindent#1}\par\medskip}
\def\datestamp{\par{\small\noindent\today}\par\vglue\bigskipamount}
\let\eref=\eqref
\let\eqref=\undefined
\allowdisplaybreaks

\def\bigo{\mathcal O}
\def\eqref#1{(\ref{#1})}
\makeatletter
\let\tru@int=\int
\def\int{\mathop{\textstyle\tru@int}\limits}
\def\overl@ss#1#2{\vcenter{\offinterlineskip
        \ialign{$\m@th#1\hfil##\hfil$\crcr#2\crcr<\crcr } }}
\def\overgr@at#1#2{\vcenter{\offinterlineskip
        \ialign{$\m@th#1\hfil##\hfil$\crcr#2\crcr>\crcr } }}
\def\gl{\mathrel{\mathpalette\overl@ss>}}
\def\lg{\mathrel{\mathpalette\overgr@at<}}

\def\Real{\mathbb{R}}

\def\Re{\mathop{\rm Re}\nolimits}

\def\pvint{\int\kern-0.94em-\kern0.2em}
\let\^=\hat
\let\==\bar
\let\@=\mathbf

\def\~#1{\tilde{#1}}
\makeatother

\def\be{\begin{equation}}
\def\ee{\end{equation}}
\def\bse{\begin{subequations}}
\def\ese{\end{subequations}}
\def\re{\mathrm{Re}~}
\def\im{\mathrm{Im}~}
\def\reg{\mathrm{reg}}

\advance\textheight 8mm
\advance\voffset -4mm

\def\br{}
\def\er{}

\begin{document}

\title{Boundary value problems for evolution partial differential equations with discontinuous data}
\author{Thomas Trogdon$^{1,*}$ and Gino Biondini$^2$}
\address{$^2$ State University of New York at Buffalo, Department of Mathematics, Buffalo, NY 14260\\
  $^1$ University of California, Irvine, Department of Mathematics, Irvine, CA 92697\\
$^*$ Corresponding author email: {\tt trogdon@math.uci.edu}}
\datestamp

\begingroup
\small\noindent
\textbf{Abstract}~
We characterize the behavior of the solutions of linear evolution partial differential equations on the half line 
in the presence of discontinuous initial conditions or discontinuous boundary conditions, 
as well as the behavior of the solutions in the presence of corner singularities.  The characterization focuses on an expansion in terms of computable special functions.
\endgroup
\bigskip

%%%%%%%%%%%%%%%%%%%%%%%%%%%%%%%%%%%%%%%%%%%%%%%%%%%%%%%%%%%%%%%%%%%%%%%%%%%%%%%%%%%%%%%%%%%%%
\section{Introduction}
\label{sec:intro}

Initial-boundary value problems (IBVPs) for linear and integrable nonlinear 
partial differential equations (PDEs)
have received renewed interest in recent years thanks to the development of the
so-called unified transform method (UTM), also known as the Fokas method.
The method provides a general framework to study these kinds of problems,
and has therefore allowed researchers to tackle a variety of interesting research questions
(e.g., see \cite{DTV,Fokas2008,fokaspelloni2014,heat,lkdv,interface,network} and references therein).

\begin{figure}[t!]
\vspace*{-2ex}
\centering
\includegraphics[width=.6\linewidth]{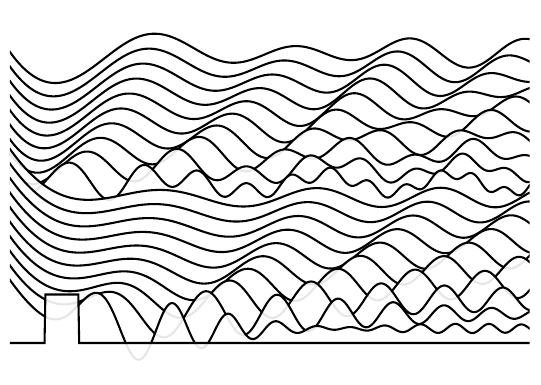}
\caption{\label{f:Airy2-simple} The solution of the Airy 2 equation \eqref{e:Airy2}  with discontinous initial and boundary data and a corner singularity.  The solution is expressed in terms of computable special functions whose asymptotics are derived in Appendix~\ref{a:SF}.  This solution is discussed in more detail in Figure~\ref{f:Airy2-sol}.}
\end{figure}

In particular, one of the topics that have been recently studied is 
that of corner singularities for IBVPs on the half line \cite{boydflyer,flyerfornberg2003b,flyerswarztrauber}.
In brief, the issue is that, on the quarter plane $(x,t)\in\Real^+\times\Real^+$, 
the limit of the PDE to the corner $(x,t)=(0,0)$ of the physical domain 
imposes an infinite number of compatibility conditions
between the initial conditions (ICs) and the boundary conditions (BCs)
[see Section~\ref{sec:UTM} for details]. 
For example, if a Dirichlet BC is given at the origin, 
the first compatibility condition is simply the requirement that the value of the IC at $x=0$ and 
that of the BC at $t=0$ are equal,
which in turn simply expresses the requirement that the solution of the IBVP 
be continuous in the limit as $(x,t)$ tends to $(0,0)$.
The higher-order compatibility conditions then arise from the repeated application of the PDE in the same limit.
Since the ICs and the BCs arise from different --- and typically independent --- domains of physics, however,
it is unlikely that they will satisfy all of these conditions.
Therefore, one could take the point of view that if one is dealing a genuine IBVP, 
one of these conditions will \textit{always} be violated.
An obvious question is then what happens when one of the compatibility conditions is violated.
Or, in other words, what is the effect on the solution of the IBVP of the violation of one among 
the infinite compatibility conditions?  See Figure~\ref{f:Airy2-simple} for an example solution where the first compatibility condition is violated and where the data is discontinuous.

Motivated by this question,
in \cite{GinoTomIVP} 
we began by considering a simpler problem. 
Namely, we studied initial value problems (IVPs) for linear evolution PDEs of the type
\begin{align}
i q_t + \omega(-i \partial _x) q = 0,
\label{e:PDE}
\end{align}
on the domain $(x,t)\in\Real\times(0,T]$, where $\omega(k)$ is a polynomial and the IC $q(x,0)$ is discontinuous.
We showed that, generally speaking, 
in the presence of dispersion and/or dissipation, the initial discontinuity is smoothed out as soon as $t\ne0$.
On the other hand, the discontinuity of the IC affects the behavior of the solution at small times.
We characterized the short-time asymptotics of the solution of the IVP
in terms of generalizations of the classical special functions, and we demonstrated 
a surprising result: Namely, that the actual solution of linear evolution PDEs with discontinuous
ICs displays all the hallmarks of the classical Gibbs phenomenon.
Explicitly: 
(i) the convergence of the solution $q(x,t)$ to the IC as $t\downarrow0$ is non-uniform
[as it should be, since $q(x,t)$ is continuous while the IC is not];
(ii) in the neighborhood of a discontinuity at $(c,0)$, the solution display high-frequency oscillations\footnote{These oscillations are characterized by a similarity solution which is obtained from the special functions.};
(iii) the oscillations are characterized by a finite ``overshoot'', which does not vanish in the limit $t\downarrow0$,
and whose value tends precisely to the Gibbs-Wilbraham contant in some appropriate limit.  This study was closely related to the work of DiFranco and McLaughlin \cite{DiFranco2005}.

In the present work we build on those results to characterize the solution of IBVPs with discontinuous data (see \cite{Lawley} for an application).
Namely, we consider the singularity propogation and smoothing properties of the linear evolution PDE 
in the domain $(x,t) \in \mathbb R^+ \times (0,T]$ with appropriate boundary data. 
Specifically, we determine a small-$x$ and small-$t$ expansion of the solution in a 
neighborhood of a discontinuity in either the boundary data or initial data.  
We also look at the solution in a neighborhood of the corner $(x,t) = (0,0)$ 
when the initial data and boundary data are not compatible.  Presumably, the methodology of Taylor \cite{Taylor2006} can be used to state that the phenomenon we describe for linear problems can be extended to certain nonlinear boundary-value problems.
Unfortunately,
unlike the case of IVPs, no general theory of well-posedness exists for IBVPs for PDEs of the form~\eref{e:PDE} with discontinuous data, and our proof of validity of the solution formula in the case of discontinuous data (Appendix~\ref{a:SolutionForm}) requires this \emph{a priori}.
Thus our treatment is necessarily limited to a few representative examples.
We emphasize, however, that:
(i) these examples describe physically relevant PDEs, 
and therefore are interesting in their own right;
(ii) since we are using the UTM, the same methodology can be applied to IBVPs for arbitrary linear evolution PDEs, if one takes well-posedness for granted. 

The outline of this work is the following:
In Section~\ref{sec:UTM} we review some relevant results about IVPs and IBVPs 
that will be used in the rest of this work.
Owing to the linearity of the PDE~\eref{e:PDE},
the solution of an IBVP with general ICs and BCs can be decomposed into the sum of 
the solution of an IBVP with the given IC and zero BCs and
the solution of an IBVP with the given BCs and zero IC. 
In Section~\ref{sec:zeroBC} we therefore characterize the solution of IBVPs with zero BCs. 
In Section~\ref{sec:zeroIC} we characterize the solution of IBVPs with zero ICs.
In Section~\ref{sec:spectraldata} we extend the results of the previous sections
to more general discontinuities.
Then, in Section~\ref{sec:examples}, we combine the results of the previous sections and 
discuss the behavior of solutions of IBVPs with corner singularities, 
i.e., the case when both ICs and BCs are non-zero but one of the compatibility conditions is violated.
%Section~\ref{sec:conclusions} concludes this work with some final remarks.

%%% Local Variables:
%%% mode: latex
%%% TeX-master: "gibbs-ibvp"
%%% End:

\section{Preliminaries}
\label{sec:UTM}

We begin by recalling some essential results about IVPs with discontinuous data from \cite{GinoTomIVP} in Section~\ref{s:IVP};
we then review the solution of  IBVPs on the half line via the \br UTM \er \cite{Fokas2008} in Section~\ref{s:UTM}.
In Sections~\ref{s:weak} and~\ref{s:compatibility} we briefly discuss weak solutions,
we present some examples of IBVPs that will be used frequently later,
and we introduce the special functions which govern the behavior of the solutions near a discontinuity.

\subsection{IVPs with discontinuous data}
\label{s:IVP}

\begin{figure}[tbp]
\kern+\smallskipamount
\centerline{\includegraphics[width=0.435\textwidth]{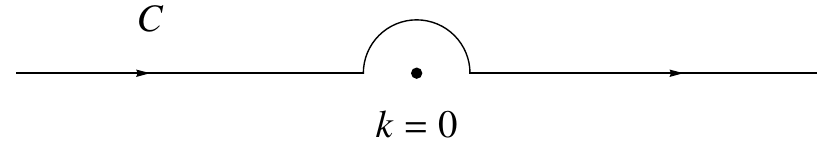}}
\kern-\smallskipamount
\caption{The integration contour $C$.}
\label{f:contour}
\end{figure}

The initial value problem  for \eref{e:PDE} with $(x,t) \in \mathbb R \times (0,T]$ 
and discontinuous ICs  was considered in~\cite{GinoTomIVP}.  The main idea  there was to consider the Fourier integral solution representation
\begin{gather}
q(x,t) = \frac{1}{2\pi} \int_{\mathbb R} e^{ikx -i\omega(k)t} \hat q_o(k) dk,
\label{e:IVPFTsoln}
\\
\hat q_o(k) = \int_{\mathbb R} e^{-ikx}q_o(x) dx, \qquad 
q(x,0) = q_o(x).
\end{gather}
Assume \br the $j$th derivative, $q_o^{(j)}$, \er has a jump discontinuity at $ x = c$.  Then $\hat q_o(k)$  can be integrated by parts to obtain
\begin{gather*}
\hat q_o(k) = e^{-ikc} \frac{[q_o^{(j)}(c)]}{(ik)^{j+1}} + \frac{F(k)}{(ik)^{j+1}},\\
[q_o^{(j)}(c)] = q_o^{(j)}(c^+) - q_o^{(j)}(c^-), \quad F(k) = \left( \int_{-\infty}^c + \int_{c}^\infty \right) q_o^{(j+1)}(x) dx.
\end{gather*}
Correspondingly,
\begin{gather}
q(x,t) = [q_o^{(j)}(c)] I_{\omega,j}(x-c,t) + \frac{1}{2\pi}\int_{C} e^{ikx - i \omega(k)t} \frac{F(k)}{(ik)^{j+1}} dk,
\label{e:IVPasympsoln} 
\\
I_{\omega,j}(x,t) = \frac{1}{2\pi} \int_C e^{ikx-i\omega(k)t} \frac{dk}{(ik)^{j+1}},
\end{gather}
where $C$ is shown in Figure~\ref{f:contour}.

The behavior of the solution formula  \eqref{e:IVPasympsoln} can then be  analyzed both near $(x,t) = (c,0)$ and near $(x,t) = (s,0), ~s \neq c$.  The function $I_{\omega,j}$   can be  examined with both the method of steepest descent and a   suitable  numerical method.  
The   second term in the right-hand side of~\eqref{e:IVPasympsoln} can be  estimated with the H\"older inequality showing that a Taylor expansion of $e^{ikx-i\omega(k)t}$ near $k = 0$ and term-by-term integration of the first $j+1$ terms produces the correct expansion (see Appendix~\ref{a:Residual}).

In this work we are concerned with the generalization of the above results to IBVPs.
The Unified Transform Method of Fokas \cite{Fokas2008} naturally   lends itself to  the above type of analysis for IBVPs, because it produces an integral representation of the solution 
  in Ehrenpreis form, similar to \eqref{e:IVPFTsoln}.

\subsection{The unified transform method for IBVPs}
\label{s:UTM}

In this section we review the Unified Transform Method (UTM) as described in \cite{Fokas2008} (see also \cite{DTV}).  The power of the method, like the Fourier transform method for pure IVPs, is that it \br gives \er an algorithmic way to produce 
an explicit integral representation of the solution, \br in Ehrenpreis form, \er of a linear, constant-coefficient IBVP on the half-line $\mathbb R^+$.

Broadly   speaking, we consider the following IBVP:
\begin{align}\label{generalPDE}
\begin{split}
iq_t + \omega(-i\partial_x)q &= 0,\quad x > 0,~t\in (0,T],\\
q(\cdot,0) &= q_o,\\
\partial_x^j q(0,\cdot) &= g_j, ~ j = 0, \ldots, N(n)-1,\\
N(n) &= \begin{cases}
n/2 & \text{$n$ even}, \\
(n+1)/2 & \text{$n$ odd and $\omega_n > 0$},\\
(n-1)/2 & \text{$n$ odd and $\omega_n< 0$},\\
\end{cases}\\
\omega(k) &= \omega_n k^n + \mathcal O(k^{n-1}).
\end{split}
\end{align}
Here $\omega(k)$ is a polynomial of degree $n$,  called the dispersion relation of the PDE.
Note that we consider the so-called canonical  IBVP, in which  the first $N(n)$ derivatives  are specified on the boundary.  
To ensure that solutions do not grow too rapidly in time, 
we impose that the imaginary part of $\omega(k)$ is bounded above.   
  In particular, in all the examples that will be discussed, $\omega(k)$ will be  real valued.  
We define the following regions in the complex $k$ plane
\begin{align*}
D = \lbrace k : \mathrm{Im}(\omega(k)) \geq 0 \rbrace, \quad D^+ = D\cap \mathbb{C}^+.
\end{align*}
 Throughout, we will 
use $L^2(I)$ to denote the space of square-integrable function on the domain $I$ and $H^k(I)$ to denote the space of functions $f$ such that $f^{(j)}$ exists a.e.~and is in $L^2(I)$ for $j = 0,1,\ldots,k$.  

Following \cite{Fokas2008,FokasSung},   one can show that,  if 
\begin{itemize}
\item $q_o \in H^{\tilde n}(\mathbb R)$, $\tilde n ={\lceil n/2 \rceil}$,
\item $g_j \in H^{1/2 + (2 \tilde n - 2j -1)/(2n)}(0,T)$ for $0 \leq j \leq N(n) -1$, and 
\item $\partial_x^j q_o(0) = g_j(0)$ for $0 \leq j \leq N(n)-1$,
\end{itemize}
then the solution of this initial-boundary-value problem is given by
\begin{align}\label{e:form}
q(x,t) = \frac{1}{2\pi}\int_{\mathbb R} e^{ikx-i\omega(k)t}\hat{q}_0(k)dk - \frac{1}{2\pi}\int_{\partial{ D^+}}\left( e^{ikx-i\omega(k)t}\sum_{j=0}^{n-1}c_j(k)\tilde{g}_j(-\omega(k),  T ) \right)dk.
\end{align}
where
\begin{align}
\hat q_o(k) = \int_{0}^\infty e^{-ikx} q_o(x) dx, \quad \tilde g_j(k,t) = \int_{0}^t e^{-iks} \partial^j_x q(0,s) ds.
\label{e:FTdefs}
\end{align}
  Hereafter, the caret 
``$\,\,\hat{\phantom{.}}\,\,$'' will refer  to the half-line Fourier transform unless specified otherwise.  
  In~\eqref{e:form}, the coefficients $c_j(k)$ are defined  by the relation
\begin{align}
i\left.\left( \frac{\omega(k) - \omega(l)}{k - l} \right)\right|_{l=-i\partial_x} = c_j(k) \partial_x^j.
\end{align}

Note that for $j > N(n)-1$, $\tilde g_j(k)$ is not specified in the statement of the problem.
Therefore, if the IBVP is well-posed, 
we expect it to be determined from the specified initial and boundary data and the PDE itself.  
This is indeed the case.  In fact,   one of the key results  of the UTM is   to show  that $\tilde g_j(k)$ can be determined purely by linear algebra.  Critical components of the theory are the so-called symmetries of the dispersion relation, 
\emph{i.e.}, the solutions $\nu(k)$ of $\omega(\nu(k)) = \omega(k)$.  
  For example, if  $\omega(k) = k^2$ then $v(k) = \pm k$ and if $\omega(k) = \pm k^3$ then $\nu(k) = k, \alpha k, \alpha^2 k$ for $\alpha = e^{2\pi i/3}$. We do not present the solution formula in any more generality.  Specifics are studied in examples.

  For our purposes, it will be convenient to  perform additional deformations to the   integration contour for the  integral along $\partial D^+$.  Let $\tilde D_i^+$, $i = 1, \ldots, N(n)$ be the connected components of $D^+$.  We   will  deform the region $\tilde D_i^+$ to a new region $D_i^+ \subset \tilde D_i^+$ such that for a given $R> 0$, $D_i^+ \cap \{ |k| < R\} = \varnothing$.  In all cases, $R$ is chosen so that all zeros of $\omega'(k)$ and $\nu(k)$ lie in the set $\{|k| < R\}$.  We display $D_i^+$ in specific examples below.

  Importantly, one can show \cite{Fokas2008} that,  for $x > 0$, $T$ in \eqref{e:form} can be replaced with $0 < t < T$ (consistently with the expectation that the solution of a true IBVP should not depend on the value of the boundary data at future times).
The replacement is not without consequences for the analysis, however.

While $\lim_{x \rightarrow 0^+} q(x,t)$ is, of course, the same in both cases, 
  the 
two formulas evaluate to give different values when computing $q(0,t)$.
This is a consequence of   the presence of  an integral in the derivation that vanishes for $x> 0$ but does not vanish for $x = 0$.  
We discuss this point more   in detail  within the context of the example \eqref{e:LS} below.  
In this work, we only study $\lim_{x \rightarrow 0^+} q(x,t)$, so this discrepancy 
is not an issue for our computations.

A similar issue is present in the evaluation of \eqref{e:form} at the point $(x,t) = (0,0)$, which is   of course  of particular interest in this work. In the case where $g_0(0) = q_o(0)$, it is apparent that neither \eqref{e:form} nor the expression obtained from \eqref{e:form} by replacing $T$ with $t = 0$ evaluates to give the correct value at the corner.  
This   issue  is discussed in more detail in the context of example \eqref{e:LS} below.  
Nevertheless, it follows from the work of Fokas and Sung \cite{FokasSung} that $\lim_{(x,t) \rightarrow (0,0)} q(x,t) = g_0(0) = q_o(0)$.  This fact also follows from our calculations.  

  The above discussion should highlight  the fact that evaluation of the solution formula near the boundary $x = 0$ and   in particular  near the corner $(x,t) = (0,0)$   of the physical domain  is indeed a non-trival matter.

\subsection{Weak solutions}
\label{s:weak}

While the Sobolev assumptions above on the initial-boundary data provide sufficient conditions 
for the representation of the solution, these assumptions must be relaxed for the purposes of the present work,
since our aim is to characterize the solution of IBVPs when either the ICs or the BCs are not differentiable.

\begin{definition}
A function $q(x,t)$  is a weak solution of \eqref{e:PDE} in an open region $\Omega$ if
\begin{align}\label{e:weak-form}
L_\omega[q,\phi] = \int_\Omega q(x,t)(-i \partial_t\phi(x,t) - \omega(i \partial_x) \phi(x,t)) dxdt = 0,
\end{align}
for all $\phi \in C^\infty_c(\Omega)$.
\end{definition}

\noindent We borrow the relaxed notion of solution of the IBVP from \cite{HolmerKdV}:
\begin{definition}
\label{d:weaksolution}
A function $q(x,t)$ is said to be an $L^2$ solution of the boundary value problem \eref{generalPDE} if 
\begin{itemize}
\item $q$ is a weak solution for $\Omega = \mathbb R^+ \times [0,T]$,
\item $q \in C^0([0,T];L^2(\mathbb R^+))$ and $q(\cdot,0) = q_o$ a.e.,
\item $\partial^j_xq \in C^0(\mathbb R^+; H^{1/2-j/n-1/(2n)}(0,T))$ and $\partial^j_xq(0,\cdot) = g_j$ a.e. for $j = 0, \ldots N(n)-1$.
%
%\marginpar{ [I don't\\know what\\``This'' refers to. - G]}
%
\end{itemize}
The \br conditions in this definition are \er obtained by setting $\tilde n = 0$.
\end{definition}

From the work of Holmer (see \cite{HolmerKdV} and \cite{HolmerNLS}) it can be inferred that when $\omega(k) = \pm k^3, \pm k^2$ the $L^2$ solutions exist and are unique.  
We are not aware of a reference that estabilishes   a similar result for more general dispersion relations, 
but we will nonetheless assume such a result to be valid.  

Two important aspects of Definition~\ref{d:weaksolution} are that (i) no compatibility conditions are required at $(x,t) = (0,0)$, 
and (ii) $H^{1/2-j/n-1/(2n)}(0,T)$ is a space that contains discontinuous functions for all $j \geq 0$.  
Another gap in the literature is that 
a set of necessary or sufficient conditions in order for \eref{e:form} to be the solution formula are, to our knowledge, not known.  
We will justify \eref{e:form} for a specific class of data that has discontinuities in Appendix~\ref{a:SolutionForm}.  
Specifically:

\begin{assume}\label{Assume:Data-0}
  The following conditions will be used in the analysis that follows:
\begin{itemize}
\item $q_o \in \mathbb L^2(\mathbb R^+) \cap L^1(\mathbb R^+, (1+|x|)^\ell)$ for some $\ell \geq 0$,
\item there exist   a sequence  $0 =x_0 < x_1 < \cdots < x_M < x_{M+1} = \infty$ such that $q_o \in H^{N(n)}( (x_i,x_{i+1}) )$ 
%for all $i =1, \ldots, M$
and $q_o(x_i^+) \neq q_o(x_i^-)$ for all $i= 1,\dots,M$, 
\item there exist   a sequence  $0 = t_0 < t_1 < \cdots < t_K < t_{K+1} = T$ such that $g_j \in H^{N(n)-j}(( t_i,t_{i+1}))$ for $i =1, \ldots, K$.
\end{itemize}
Note that $g_j$ may or may not be discontinuous at each $t_i$.
\end{assume}

Our results on sufficient conditions for \eqref{e:form} to produce the solution formula are not complete.  We consider the full development of this topic important but beyond the scope of this paper.

\subsection{Compatibility conditions}
\label{s:compatibility}

In this section we discuss the conditions required to ensure that no singularity \br is \er present at the corner $(x,t) = (0,0)$.  
The first $N(n)$ conditions are simply given by
\begin{align*}
q_o^{(j)}(0) = g_j(0), \quad j = 0, \ldots, N(n) -1.
\end{align*}
Higher-order conditions are found by enforcing that the differential equation holds at the corner:
\begin{align*}
i g_j^{(\ell)}(0) + \omega(-i\partial_x)^{\ell} q_o^{(j)}(0) = 0, \quad \ell =1,2,\ldots.
\end{align*}
We refer to the index $j+n \ell$ as the \textit{order} of the compatibility condition.  
Note that because $N(n)-1 < n$, there is not a compatibility condition at every order.  
Still, if $m$ is an integer we say that the compatibility conditions hold up to order $m$ if they hold for every choice of $j$ and $\ell$ such that $j+n\ell \leq m$.

\subsection{Examples}
\label{s:examples}

In the rest of this work we will   illustrate  our results by discussing several examples of physically relevant IBVPs.
Therefore, we recall, for convenience, the solution formulae for these IBVPs, 
as obtained with the unified transform method.
We refer the reader to Refs.~\cite{Fokas2008,fokaspelloni2014} for all details.

\subsubsection{Linear Schr\"odinger}

\iffalse
\begin{figure}[tbp]
\centering
\includegraphics[width=.6\linewidth]{figs/LS-Regions.pdf}
\caption{\label{Figure:LS-Regions} A schematic for $\omega(k) = k^2$.  The shaded region  is $D$.}
\end{figure}
\fi

Consider the IBVP
\bse
\label{e:LS}
\begin{gather}
iq_t + q_{xx} = 0, \quad  x\geq 0,~t\in(0,T],\\
q(\cdot,0) = q_o,\qquad
q(0,\cdot) = g_0.
\end{gather}
\ese
The dispersion relation is $\omega(k) = k^2$, and the solution formula for the IBVP is given by (replacing $T$ with $t$) in \eqref{e:form})
\begin{align*}
q(x,t) = \frac{1}{2 \pi} \int_{\mathbb R} e^{ikx-i\omega(k)t} \hat q_o(k) dk + \frac{1}{2\pi} \int_{\partial D^+} e^{ikx-i\omega(k)t} [2k \tilde g_0(-\omega(k),t)-\hat q_o(-k)]dk.
\end{align*}
See Figure~\ref{Figure:LS-Regions} for $D^+$ and $D_1^+$.

\begin{figure}[t!]
\begin{center}
\begin{minipage}[t!]{0.45\textwidth}%
\includegraphics[width=\linewidth]{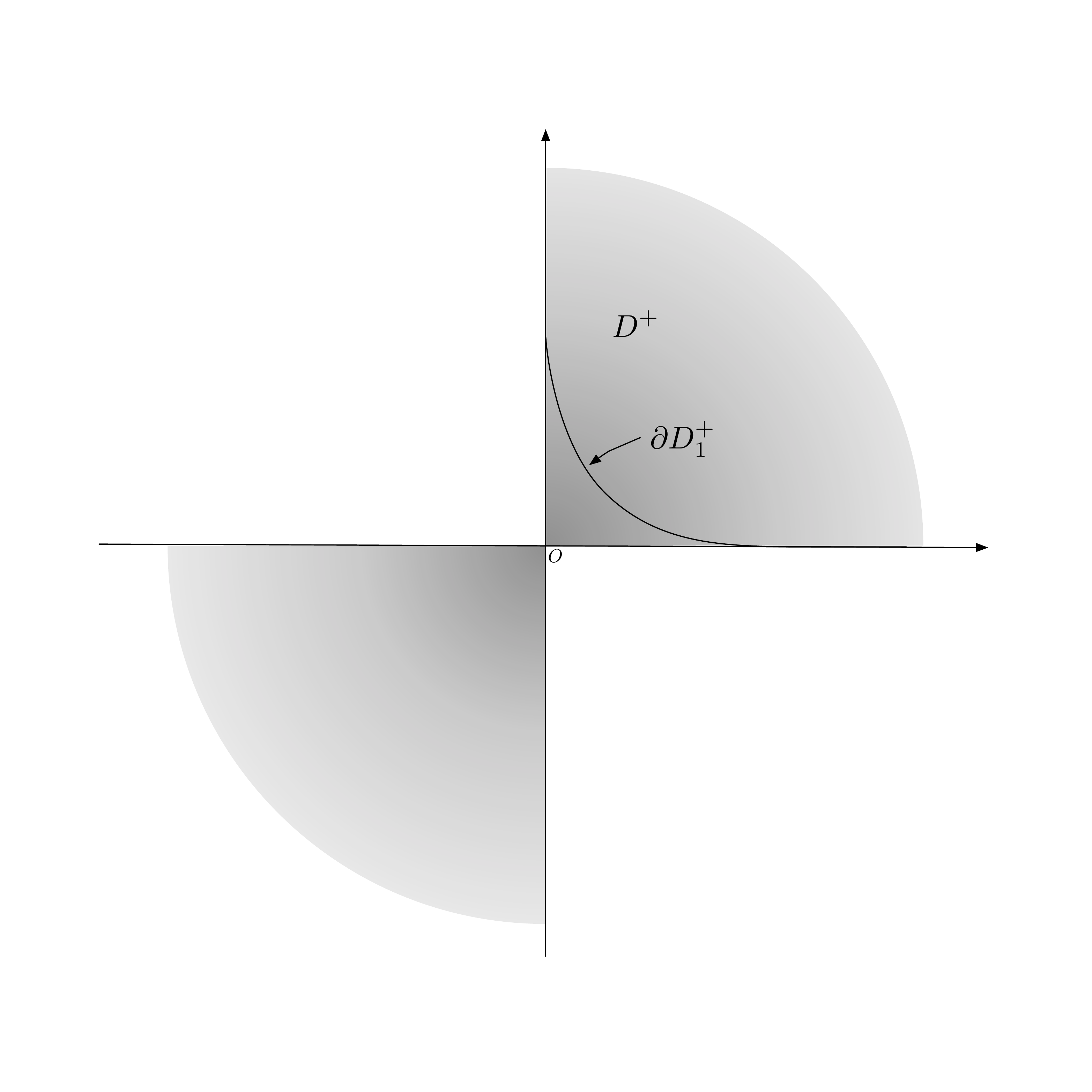}
\caption{The region $D$ (shaded) in the complex $k$-plane for the linear Schr\"odinger equation~\eqref{e:LS}, corresponding to $\omega(k) = k^2$.
The modified contour $\partial D_1^+$ is also shown.
}
\label{Figure:LS-Regions} 
\end{minipage}~~~~~
\begin{minipage}[t!]{0.45\textwidth}%
\includegraphics[width=\linewidth]{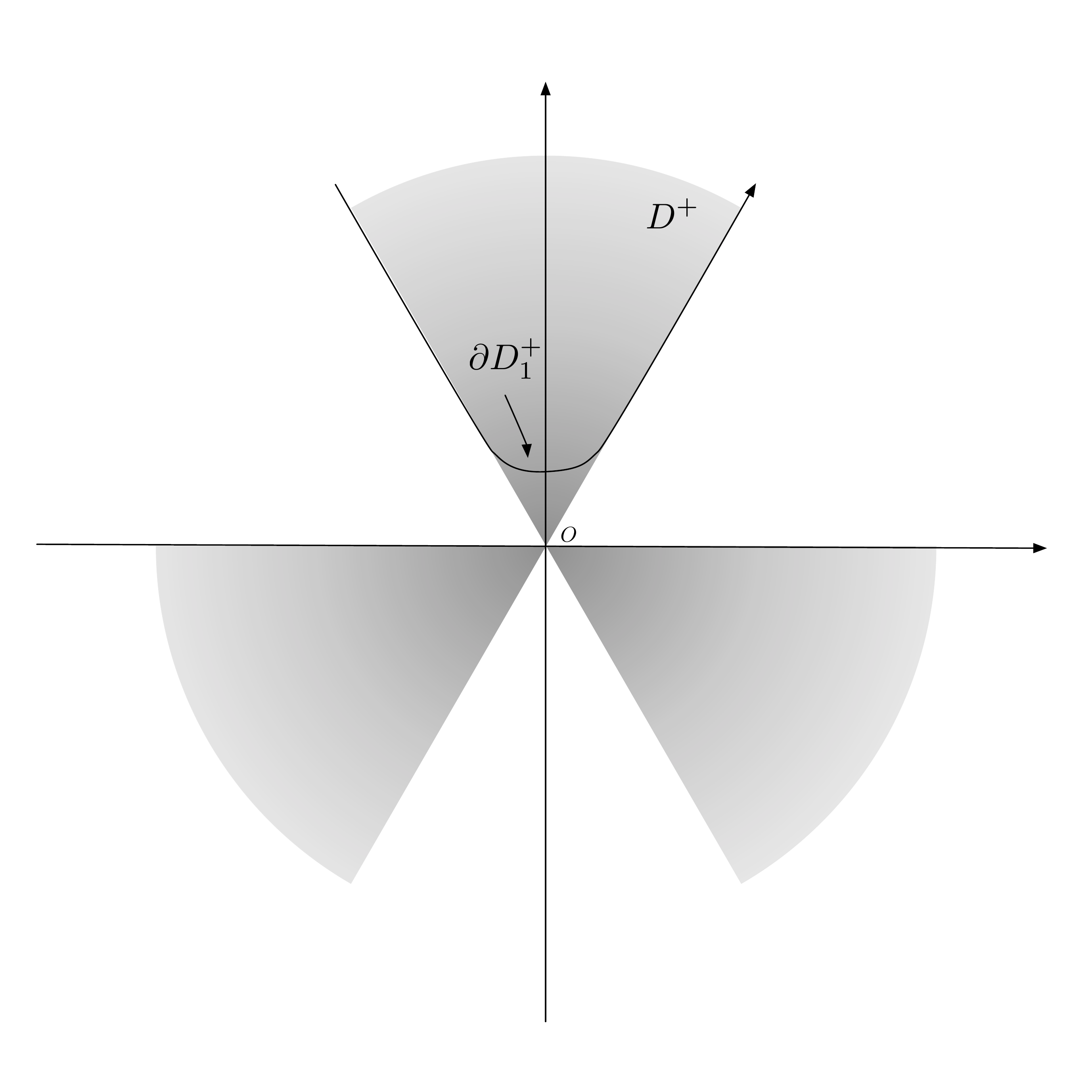}
\caption{  Same as Fig.~\ref{Figure:LS-Regions}, but  for the Airy~1 equation~\eqref{e:Airy1}, corresponding to $\omega(k) = -k^3$.\protect\newline~\protect\newline~}
\label{Figure:Airy1-Regions} 
\end{minipage}\\
\begin{minipage}[t!]{0.45\textwidth}%
\includegraphics[width=\linewidth]{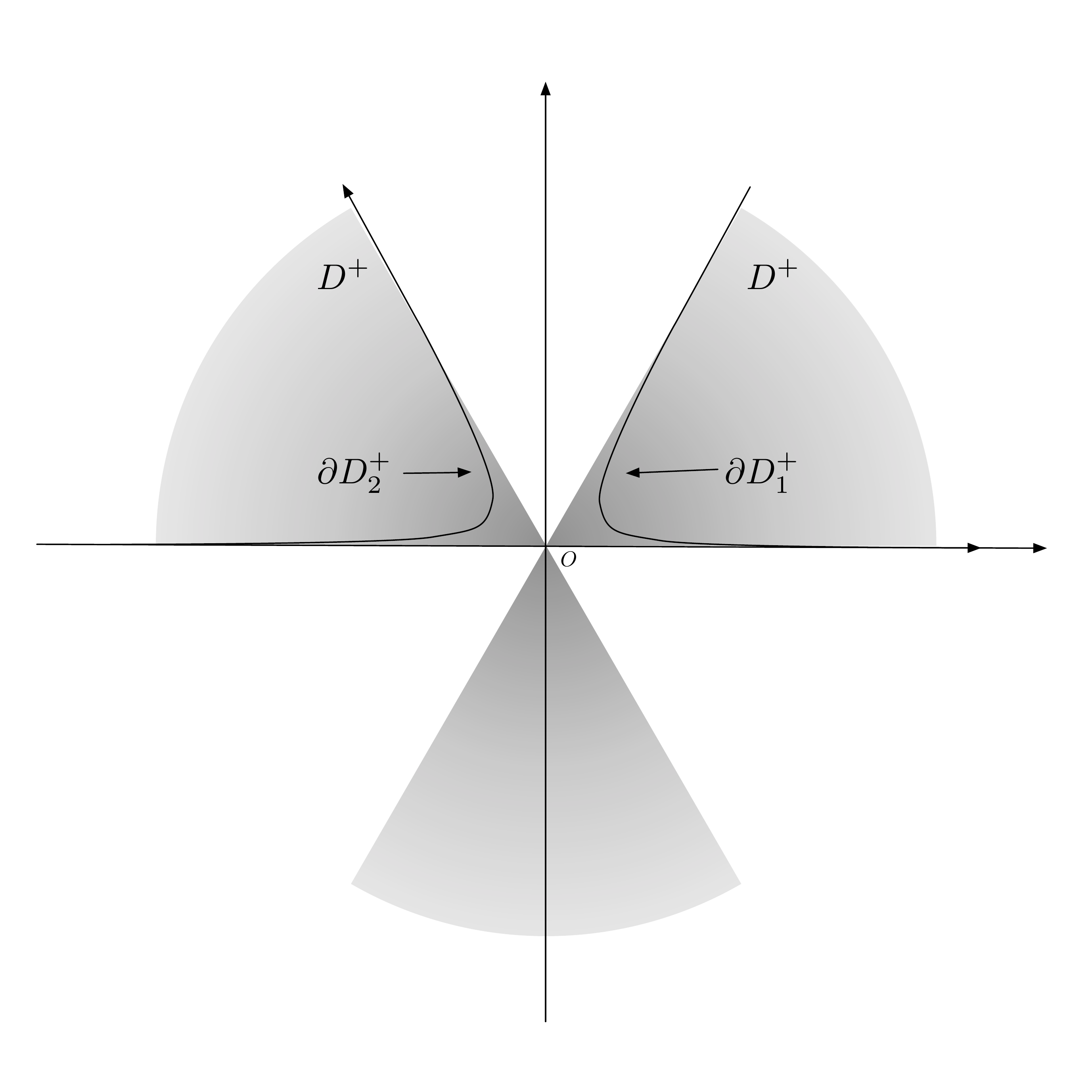}
\caption{  Same as Fig.~\ref{Figure:LS-Regions}, but  for the Airy~2 equation~\eqref{e:Airy2}, corresponding to $\omega(k) = k^3$.
\protect\newline~\protect\newline~}
\label{Figure:Airy2-Regions}
\end{minipage}
\end{center}
\end{figure}

For this specific example, we discuss the evaluation of $q(x,t)$ at $x = 0$ and at $(x,t) = (0,0)$ in detail,
  in order to illustrate some of the issues that arise when taking the limit of the solution representation~\eqref{e:form}. 
We assume continuity of $q_o$ and $g_0$ and rapid decay of $q_o$ at infinity.  First, by contour deformations, for $t > 0$, the solution formula   can be written  as
\begin{align}\label{e:t}
  q(0,t) = \frac{1}{2 \pi} \int_{\mathbb R} e^{-i\omega(k)t} [\hat q_o(k) - \hat q_o(-k)] dk - \frac{1}{2\pi} \int_{\partial D^+} e^{-i\omega(k)t} 2k \tilde g_0(-\omega(k),t) dk.
\end{align}
Then by the change of variables $k \mapsto -k$, the first integral   can be shown to vanish  identically. For this last integral we let $s = - \omega(k)$ and find
\begin{align*}
  q(0,t) &=  \frac{1}{2\pi} \int_{-\infty}^\infty e^{ist} \tilde g_0(s,t) ds = \frac{1}{2\pi} \int_{-\infty}^\infty e^{ist} \left(\int_{0}^t e^{-i\tau s} g_0(\tau) d \tau\right) ds \\
  &= \frac{1}{2\pi} \int_{-\infty}^\infty e^{ist} \left(\int_{-\infty}^\infty e^{-i\tau s} g_0(\tau) \chi_{[0,t]}(\tau) d \tau\right) ds = \half g_0(t).
\end{align*}
Here we use $g_0(\tau) \chi_{[0,t]}(\tau) = 0$ for $\tau \not\in [0,t]$ and $\half g_0(t)$ is the average value of the left and right limits of this function at $\tau = t$.  If $T$ is used in \eqref{e:form} and $t < T$, we have
\begin{align}\label{e:T}
  q(0,t) = \frac{1}{2\pi} \int_{-\infty}^\infty e^{ist} \left(\int_{-\infty}^\infty e^{-i\tau s} g_0(\tau) \chi_{[0,T]} d \tau\right) ds = g_0(t),
\end{align}
because $g_0(\tau) \chi_{[0,T]}(\tau)$ is continuous at $\tau = t$. Now, by similar arguments, if $t = 0$ we get \br zero  for \eqref{e:T} and the first integral in \eqref{e:t}\er. Nevertheless, the limit to the boundary of the domain from the interior produces the correct values.

\subsubsection{Airy 1}

\iffalse
\begin{figure}[tbp]
\centering
\includegraphics[width=.6\linewidth]{figs/Airy1-Regions.pdf}
\caption{\label{Figure:Airy1-Regions} A schematic for $\omega(k) = -k^3$.  The shaded region  is $D$.}
\end{figure}
\fi

Consider the IBVP
\bse
\label{e:Airy1}
\begin{gather}
q_t + q_{xxx} = 0, \quad  x\geq 0,~t\in(0,T],\\
q(\cdot,0) = q_o,\qquad
q(0,\cdot) = g_0.
\end{gather}
\ese
The dispersion relation is $\omega(k) = -k^3$, and the solution of the IBVP is given by
\begin{align}
q(x,t) &=\frac{1}{2\pi}\int_{-\infty}^\infty e^{ikx-i\omega(k)t}\hat{q}_0(k)dk - \frac{1}{2\pi}\int_{\partial D^+}3k^2e^{ikx-i\omega(k)t}\tilde{g}_0(-\omega(k),t) dk \nonumber \\
&+ \frac{1}{2\pi}\int_{\partial D^+}e^{ikx-i\omega(k)t} \left[\alpha\hat{q}_0(\alpha k) + \alpha^2\hat{q}_0(\alpha^2 k) \right]dk.
\end{align}
See Figure~\ref{Figure:Airy1-Regions} for $D^+$ and $D_1^+$.

\subsubsection{Airy 2}

\iffalse
\begin{figure}[tbp]
\centering
\includegraphics[width=.6\linewidth]{figs/Airy2-Regions.pdf}
\caption{\label{Figure:Airy2-Regions} A schematic for $\omega(k) = k^3$.  The shaded region  is $D$.}
\end{figure}
\fi

Consider the IBVP
\bse
\label{e:Airy2}
\begin{gather}
q_t - q_{xxx} = 0, \quad  x > 0,~t\in(0,T],\\
q(\cdot,0) = q_o,\qquad
q(0,\cdot) = g_0,\qquad
q_x(0,\cdot)  = g_1.
\end{gather}
\ese
Note that two BCs need to be assigned at $x=0$, unlike the previous example.  The dispersion relation is $\omega(k) = k^3$,
and the integral representation for the solution of the IBVP is
\begin{align}\label{e:Airy2-sol}
q(x,t) = \frac{1}{2\pi}\int_{\mathbb R} e^{ikx-i\omega(k)t}\hat{q}_0(k)dk - \frac{1}{2\pi}\int_{\partial D_1^+} e^{ikx-i\omega(k)t}\tilde{g}(k,t)dk - \frac{1}{2\pi}\int_{\partial D_2^+} e^{ikx-i\omega(k)t}\tilde{g}(k,t)dk,
\end{align}
where
\bse
\begin{align}
\tilde{g}(k,t) &= \hat{q}_0(\alpha k) + (\alpha^2-1)k^2\tilde{g}_0(-\omega(k),t) - i(\alpha-1) k\tilde{g}_1(-\omega(k),t),& ~~k\in \partial D_2^+, \\
\tilde{g}(k,t) &= \hat{q}_0(\alpha^2 k) + (\alpha-1) k^2\tilde{g}_0(-\omega(k),t) - i(\alpha^2-1) k\tilde{g}_1(-\omega(k),t),& ~~k\in \partial D_1^+.
\end{align}
\ese
See Figure~\ref{Figure:LS-Regions} for $D^+$, $D_1^+$ and $D_2^+$.

\subsection{Special functions}
\label{s:special}

In the following we will make extensive use of the functions

\begin{align}
I_{\omega,m,j}(x,t) = \frac{1}{2\pi} \int_{\partial D_j^+} e^{ikx-i\omega(k)t} \frac{dk}{(ik)^{m+1}}.
\end{align}
Also, when taking the sum over all contours
we use the modified notation
\begin{align}
I_{\omega,m}(x,t) = \frac{1}{2\pi} \int_{C} e^{ikx-i\omega(k)t} \frac{dk}{(ik)^{m+1}} = \sum_{j=1}^{N(n)} I_{\omega,m,j}(x,t).
\end{align}
The properties of these functions are discussed in Appendix~\ref{a:SF}.

%%% Local Variables:
%%% mode: latex
%%% TeX-master: "gibbs-ibvp"
%%% End:

\section{IBVP with zero boundary data}
\label{sec:zeroBC}

By Lemma~\ref{Lemma:Valid}, we know that 
the solution formula \eqref{e:form} holds for piecewise smooth data without any compatibility conditions imposed at $x = 0$, $t =0$.  We begin with assuming zero boundary data and then we relax our assumptions systematically.  We perform this analysis on a case-by-case basis and then generalize our results.  There are four relevant components of the analysis of this solution formula:
\begin{enumerate}
\item the behavior of $q$ near $x = 0$ for $t > 0$,
\item the behavior of $q$ near $(x,t) = (0,0)$,
\item the behavior of $q$ near $(x,t) = (c,0)$ when $c$ is a discontinuity of $q_o$, and
\item the behavior of $q$ near $(x,t) = (s,0)$ when $q_o$ is continuous at $s$.
\end{enumerate}

\subsection{Linear Schr\"odinger}

With zero Dirichlet BCs, the solution of \eref{e:LS} is given by (recall $\omega(k) = k^2$)
\begin{align*}
q(x,t) = \frac{1}{2 \pi} \int_{\mathbb R} e^{ikx-i\omega(k)t} \hat q_o(k) dk - \frac{1}{2\pi} \int_{\partial D^+} e^{ikx-i\omega(k)t} \hat q_o(-k)dk.
\end{align*}
In this simple case, the solution can be found by a straightforward application of the method of images.  Also, the integral on $\partial D^+$ can be deformed back to the real axis.  However, below we will encounter situations where this deformation is not possible, 
so we will treat this case by keeping the second	 integral on $\partial D^+$.

\subsubsection{Short-time behavior}

We consider Assumption~\ref{Assume:Data-0} with $g_0 \equiv 0$ (i.e., zero BC).  
 We begin by studying the case  $M =0$, \emph{i.e.}, the IC has no discontinuities in $\mathbb R^+$. 
 On the other hand, if $q(0)\ne0$, the compatibility condition at $(x,t)=(0,0)$ is not satisfied. 
As discussed in the introduction, we integrate  the first of~\eqref{e:FTdefs}  by parts,  to obtain 
\vspace*{-1ex}
\begin{align}
\hat q_o(k) &= \frac{q_o(0)}{ik} + \frac{F_0(k)}{ik}, \quad 
F_0(k) = \int_{0}^\infty e^{-ikx} q_o'(x) dx.
\label{e:LS_F0def}
\end{align}
After a contour deformation, we are then left with 
\begin{align}\label{e:ibp-LS}
q(x,t) = 2 q_o(0) I_{\omega,0,1}(x,t) + \frac{1}{2 \pi} \int_{C} e^{ikx-i\omega(k)t} \frac{F_0(k)}{ik} dk + \frac{1}{2 \pi} \int_{\partial D^+_1} e^{ikx-i\omega(k)t} \frac{F_0(-k)}{ik} dk.
\end{align}
We now appeal to Lemmas~\ref{Lemma:Expansion} and \ref{Lemma:Sufficient} to derive an expansion about $(s,0)$:
\begin{align}
q(x,t) = 2 q_o(0) I_{\omega,0,1}(x,t) + \frac{1}{2 \pi} \int_{C} e^{iks} \frac{F_0(k)}{ik} dk + \frac{1}{2 \pi} \int_{\partial D^+_1} e^{iks} \frac{F_0(-k)}{ik} dk + \bigo( |x-s|^{1/2} + t^{1/4}).
\label{e:LSexpn1}
\end{align}

\begin{remark}
The expansion~\eqref{e:LSexpn1} can be interpreted by noting that,  in a neighborhood of $(s,0)$, the difference
\vspace*{-1ex}
\begin{align*}
q(x,t) - 2 q_o(0) I_{\omega,0,1}(x,t)
\end{align*}
has an expansion in terms of functions depending only on $s$, up to the error terms, and hence the leading-order $x$-dependence  of $q(x,t)$  is captured by $2 q_o(0) I_{\omega,0,1}(x,t)$.
\end{remark}

It follows  from~\eqref{e:LS_F0def}  that $F_0(k)$ is analytic and decays in the lower-half plane, so that $F_0(-k)$ has the same properties in the upper-half plane.  This implies  that for all $s>0$, the third term in the right-hand side of~\eqref{e:LSexpn1} vanishes identically: 
\begin{align*}
\frac{1}{2 \pi} \int_{\partial D^+_1} e^{iks} \frac{F_0(-k)}{ik} dk = 0, \quad s > 0.
\end{align*}
Furthermore, if $s \neq 0$ one can use Theorem~\ref{Thm:Kernel} to work out the behavior of $I_{\omega,0,1}$ to obtain, 
noting that its error term is $\mathcal O(t^{1/2})$:
\begin{align*}
q(x,t) = \frac{1}{2 \pi} \int_{C} e^{iks} \frac{F_0(k)}{ik} dk + \bigo( |x-s|^{1/2} + t^{1/4}),\quad s > 0\,.
\end{align*}
As expected, this is the same behavior as for the IVP (see \cite{GinoTomIVP}): We recover the initial condition.

Next, we  consider the case  $q_o(0) = 0$ and $M = 1$,  
 implying that the first compatibility condition at $(x,t)=(0,0)$ is satisfied, but the IC is discontinuous at $x=x_1$.  Again, integration by parts produces
\begin{align*}
q(x,t) &=  [q_o(x_1)] I_{\omega,0,1}(x-x_1,t) + [q_o(x_1)] I_{\omega,0,1}(x+x_1,t)  \\
&+ \frac{1}{2 \pi} \int_{C} e^{ikx-i\omega(k)t} \frac{F_0(k)}{ik} dk + \frac{1}{2 \pi} \int_{\partial D^+_1} e^{ikx-i\omega(k)t} \frac{F_0(-k)}{ik} dk,
\end{align*}
and Lemmas~\ref{Lemma:Expansion} and \ref{Lemma:Sufficient} produce an expansion 
\begin{align*}
q(x,t) &=  [q_o(x_1)] I_{\omega,0,1}(x-x_1,t) + \frac{1}{2 \pi} \int_{C} e^{iks} \frac{F_0(k)}{ik} dk + \frac{1}{2 \pi} \int_{\partial D^+_1} e^{iks} \frac{F_0(-k)}{ik} dk + \bigo(|x-s|^{1/2} + t^{1/4})
\end{align*}
for all $s\geq0$.
Here Theorem~\ref{Thm:Kernel} was used to discard $I_{\omega,0,1}(x+x_1,t)$ (its error term is smaller, $\mathcal O(t^{1/4})$).  Continuing, for $s \neq x_1$, $s \neq 0$  we also have 
\begin{align*}
q(x,t) &=  -[q_o(x_1)] \chi_{(-\infty,0)}(x-x_1) + \frac{1}{2 \pi} \int_{C} e^{iks} \frac{F_0(k)}{ik} dk + \bigo(|x-s|^{1/2} + t^{1/4}),
\end{align*}
 from which  additional considerations (cf.\ \cite{GinoTomIVP}) yield 
\begin{align*}
q(x,t) &= q_o(s) + \bigo(|x-s|^{1/2} + t^{1/4}),
\end{align*}
as is expected.

%\marginpar{[What's the second way? I assumed it's the cut-off - G]}
%
\br The general case can be explained in the following way. First, we point out that we not employ \er  integration by parts on each interval of differentiability of $q_o$.   In the IVP,  the difficulty in using this  approach  is that  one needs to implicitly assume  analyticity of $\hat q_o(k)$ throughout (in order  to deform  the integration contour  to $C$) and the  requirements needed to ensure analyticity  are too restrictive. 
In the IBVP,  on the other hand, we automatically  have analyticity  for $\hat q_o(k)$ in the lower-half plane,  so this  requirement  is not an issue.
But  in order to keep our treatment consistent with that for  the IVP \br (as in \cite{GinoTomIVP}) \er we use cut-off functions.  Let $\phi_\epsilon(x)$ be supported on $[-\epsilon,\epsilon]$, equal to unity for $x \in [-\epsilon/2,\epsilon/2]$ and interpolate monotonically and infinitely smoothly between 0 and 1 on $[-\epsilon,-\epsilon/2)$ and $(\epsilon/2,\epsilon]$.  Examples of such functions are well-known \cite{bleisteinhandlesman} (see also \cite{GinoTomIVP}).  We decompose the initial condition as follows
\begin{align*}
q_o(x) = \sum_{m=0}^M \underbrace{q_o(x) \phi_\epsilon(x-x_m)}_{q_{o,m}(x)} +  \underbrace{q_o(x) \left( 1 - \sum_{m=1}^M \phi_{\epsilon}(x-x_m) \right)}_{q_{o,\reg}(x)},
\end{align*}
with $\epsilon < \min_{m=1,\dots,M-1} |x_m-x_{m+1}|/2$.

 Each of  the Fourier transforms $\hat q_{o,m}$ is analytic near $k = 0$, so a deformation  of the integration contour  to $C$  for each of them  is  now  justified.  The results of this section produce asymptotics of the solutions $q_{j}(x,t)$  obtained with each of  these initial conditions.  It remains to understand the behavior of $q_{\reg}(x,t)$.  
 If one extends $q_{o,\reg}$  to be zero for $x < 0$,
one has $q_{o,\reg} \in H^1(\mathbb R)$,
which implies
\begin{align*}
q_\reg(x,t) -q_{o,\reg}(s)  = \bigo(|x-s|^{1/2} + t^{1/4}).
\end{align*}
Combining everything  \br we then have  
\begin{align*}
q(x,t) &= \frac{1}{2 \pi} \int_{C} e^{iks} \frac{F_0(k)}{ik} dk + \frac{1}{2 \pi} \int_{\partial D^+_1} e^{iks} \frac{F_0(-k)}{ik} dk  + \sum_{ m=1}^{ M} [q_o(x_m)] I_{\omega,0,1}(x-x_m,t)\\
&+ 2 q_o(0)I_{\omega,0,1}(x,t) +\bigo(|x-s|^{1/2} + t^{1/4}), \quad F_0(k) = \int_0^\infty e^{iks} \frac{d}{ds}[ q_o(s) - q_{o,\reg}(s)] ds,
\end{align*}
where the differentiation in this last line occurs on each interval of differentiability. \er  Note that the integral on $\partial D_1^+$ vanishes when $s > 0$.

\subsubsection{Boundary behavior}

It is straightforward to check that $q(0,t) = 0$ for $t > 0$. If $\ell$ is sufficiently large in the sense of Theorem~\ref{Thm:Reg}, then Taylor's theorem implies $q(x,t) = \bigo(x)$ for $t \geq \delta > 0$.  If only $L^2$ assumptions are made, then Lemmas~\ref{Lemma:Expansion} and \ref{Lemma:Sufficient} with the above expansion produce $|q(x,t)| \leq C |x|^{1/2}$ where $C$ depends on $\|q_o'\|_{L^2(\mathbb R^+)}$ and $[q_o(x_i)] I_{\omega,0,1}(x-x_i,t)$.  This derivative is taken to be defined piecewise on its intervals of differentiability.

\subsection{Airy 1}

With zero boundary data, we consider the solution of \eref{e:Airy1} (Assumption~\ref{Assume:Data-0} with $g_0 \equiv 0$)
\begin{align}
q(x,t) &=\frac{1}{2\pi}\int_{\mathbb R} e^{ikx-i\omega(k)t}\hat{q}_0(k)dk + \frac{1}{2\pi}\int_{\partial D^+}e^{ikx-i\omega(k)t} \left(\alpha\hat{q}_0(\alpha k) + \alpha^2\hat{q}_0(\alpha^2 k) \right)dk.
\end{align}

\subsubsection{Short-time behavior}

We proceed as before.  First, assume the initial data is continuous (again, along with Assumption~\ref{Assume:Data-0}).  After integration by parts,  we must consider the integral ($\omega(k) = -k^3$)
\begin{align*}
q(x,t) = 3 q_o(0) I_{\omega,0,1}(x,t) + \frac{1}{2 \pi} \int_{C} e^{ikx-i\omega(k)t} F_0(k) \frac{dk}{ik} + \frac{1}{2 \pi} \int_{\partial D^+_1} e^{ikx-i\omega(k)t} (F_0( \alpha k) + F_0(\alpha^2 k)) \frac{dk}{ik}.
\end{align*}
The analysis of this expression is not much different from \eref{e:ibp-LS}.  Next, we assume $q_o(0) = 0$, $M = 1$.  We obtain
\begin{align*}
q(x,t) &= [q_o(x_1)]  \left( I_{\omega,0,1}(x-x_1,t) + I_{\omega,0,1}(x-\alpha x_1,t) + I_{\omega,0,1}(x-\alpha ^2x_1,t) \right) \\
&+ \frac{1}{2 \pi} \int_{C} e^{iks} F_0(k) \frac{dk}{ik} + \frac{1}{2 \pi} \int_{\partial D^+_1} e^{iks} (F_0( \alpha k) + F_0(\alpha^2 k)) \frac{dk}{ik} + \bigo(|x-s|^{1/2} + t^{1/6})
\end{align*}
by appealing to Lemmas~\ref{Lemma:Expansion} and \ref{Lemma:Sufficient}.  More care is required to understand $I_{\omega,0,1}(x- \alpha x_1)$.  Specifically, we look at $
e^{ik(x-\alpha x_1)}$ and $\partial D_1^+$.  For sufficiently large $k \in \partial D_1^+$, $k = \pm |k| \cos \theta + i |k| \sin \theta$ for $\theta = 2 \pi /3$. \br It \er follows that
\begin{align*}
 \br \Re i k (x -\alpha x_1) = - |k| x_1 \sin \theta + |k| x_1 \sin (\theta + \phi) \leq 0, \quad \text{for } x \geq 0, ~~ \frac{\theta}{2} \leq \phi \leq \theta. \er
\end{align*}
Jordan's Lemma can be applied to show that $I_{\omega,0,1}(x-\alpha x_1,0) = 0$ for $x \geq 0$.  We write
\begin{align*}
I_{\omega,0,1}(x-\alpha x_1,t) = \frac{1}{2 \pi} \int_{\partial D_1^{+}} ( e^{- i \omega(k) t} -1)e^{ik(x - \alpha x_1)}\frac{dk}{ik} = \bigo(t^{1/6}),
\end{align*}
by appealing to Lemma~\ref{Lemma:Expansion}.  Similar calculations hold for $I_{\omega,0,1}(x- \alpha^2 x_1,t)$.  Therefore,
\begin{align*}
q(x,t) &= [q_o(x_1)]  I_{\omega,0,1}(x-x_1,t) + \frac{1}{2 \pi} \int_{C} e^{iks} F_0(k) \frac{dk}{ik} \\
&+ \frac{1}{2 \pi} \int_{\partial D^+_1} e^{iks} (F_0( \alpha k) + F_0(\alpha^2 k)) \frac{dk}{ik} + \bigo(|x-s|^{1/2} + t^{1/6}).
\end{align*}
If $s \neq x_1$, $I_{\omega,0,1}(x-x_1,t)$ can be replaced with $-\chi_{(-\infty,0)}(x-x_1)$.  Finally, if $s > 0$ then the integral on $\partial D^+_1$ vanishes identically.  Combining everything, in the general case we have
\begin{align*}
q(x,t) &= \frac{1}{2 \pi} \int_{C} e^{iks} F_0(k) \frac{dk}{ik} + \frac{1}{2 \pi} \int_{\partial D^+_1} e^{iks} (F_0( \alpha k) + F_0(\alpha^2 k)) \frac{dk}{ik}  + \sum_{x_i} [q_o(x_i)] I_{\omega,0,1}(x-x_i,t)\\
&+ 3q_o(0) I_{\omega,0,1}(x,t)+ \bigo(|x-s|^{1/2} + t^{1/6}).
\end{align*}
Here the integral on $\partial D^+_1$ should be dropped with $s > 0$.

\begin{remark}
Again, this expansion is interpreted by noting that
\begin{align*}
q(x,t) - \sum_{x_i} [q_o(x_i)] I_{\omega,0,1}(x-x_i,t) - 3q_o(0) I_{\omega,0,1}(x,t)
\end{align*}
has an expansion in a neighborhood of $(s,0)$ in terms of smoother functions depending only on $s$, \br where $s$ is fixed,\er  up to the error terms.
\end{remark}

\subsubsection{Boundary behavior}
Finally, for $x=0$, $t \geq \delta > 0$, we know the solution is smooth from Theorem~\ref{Thm:Reg} and $q(0,t) = 0$ so that we find $q(x,t) = \bigo(x)$ from Taylor's theorem.  Again, $q(x,t) = \bigo(x^{1/2})$ follows if only $L^2$ assumptions are made on the initial data and its derivative.

\subsection{Airy 2}

Recall that the solution to \eref{e:Airy2} is given by \eref{e:Airy2-sol} with ($\omega(k) = k^3$)
\begin{align}
\tilde{g}(k,t) &= \hat{q}_0(\alpha k), ~~k\in D_2^+, \\
\tilde{g}(k,t) &= \hat{q}_0(\alpha^2 k), ~~k\in D_1^+.
\end{align}
when the boundary data is set to zero.

\subsubsection{Short-time behavior}

 Following the same procedure, we assume the initial data is continuous and find
\begin{align*}
q(x,t) &= q_o(0) \left(I_{\omega,0}(x,t) - \alpha^{-1} I_{\omega,0,2}(x,t) - \alpha^{-2} I_{\omega,0,1}(x,t) \right) + \frac{1}{2 \pi} \int_{C} e^{ikx-i\omega(k)t} \frac{F_0(k)}{ik} dk \\
&- \frac{\alpha^{-1}}{2 \pi} \int_{\partial D^+_2} e^{ikx-i\omega(k)t} \frac{F_0(\alpha k)}{ik} dk - \frac{\alpha^{-2}}{2 \pi} \int_{\partial D^+_1} e^{ikx-i\omega(k)t} \frac{F_0(\alpha^2 k)}{ik} dk.
\end{align*}
Then  the expansion
\begin{align*}
q(x,t) &= q_o(0) \left(I_{\omega,0}(x,t) - \alpha^{-1} I_{\omega,0,2}(x,t) - \alpha^{-2} I_{\omega,0,1}(x,t) \right) + \frac{1}{2 \pi} \int_{C} e^{iks} \frac{F_0(k)}{ik} dk \\
&- \frac{\alpha^{-1}}{2 \pi} \int_{\partial D^+_2} e^{iks}\frac{F_0(\alpha k)}{ik} dk - \frac{\alpha^{-2}}{2 \pi} \int_{\partial D^+_1}e^{iks} \frac{F_0(\alpha^2 k)}{ik} dk + \bigo(|x-s|^{1/2} + t^{1/6}).
\end{align*}
follows.  If $s > 0$ the first three terms may be removed.  Furthermore, the terms involving $F_0(\alpha k)$ and $F_0(\alpha^2 k)$ vanish identically if $s > 0$.  Now, assume $q_o(0) = 0$ and $M =1$.  We find
\begin{align*}
  q(x,t) &=  [q_o(x_1)] \left(I_{\omega,0}(x-x_1,t) - \alpha^{-1} I_{\omega,0,2}(x-\alpha x_1,t) - \alpha^{-2} I_{\omega,0,1}(x-\alpha^2 x_1,t) \right) \\
&+ \frac{1}{2 \pi} \int_{C} e^{ikx-i\omega(k)t} \frac{F_0(k)}{ik} dk - \frac{\alpha^{-1}}{2 \pi} \int_{\partial D^+_2} e^{ikx-i\omega(k)t} \frac{F_0(\alpha k)}{ik} dk - \frac{\alpha^{-2}}{2 \pi} \int_{\partial D^+_1} e^{ikx-i\omega(k)t} \frac{F_0(\alpha^2 k)}{ik} dk\\
&= [q_o(x_1)]  I_{\omega,0}(x-x_1,t)+ \frac{1}{2 \pi} \int_{C} e^{iks} \frac{F_0(k)}{ik} dk \\
&- \frac{\alpha^{-1}}{2 \pi} \int_{\partial D^+_2} e^{iks}\frac{F_0(\alpha k)}{ik} dk - \frac{\alpha^{-2}}{2 \pi} \int_{\partial D^+_1} e^{iks}\frac{F_0(\alpha^2 k)}{ik} dk + \bigo(|x-s|^{1/2} + t^{1/6}),
\end{align*}
because it can be shown that $I_{\omega,0,2}(x-\alpha x_1,t)  = I_{\omega,0,1}(x-\alpha^2 x_1,t) = \bigo(t^{1/6})$ for $x > 0$ in the same way as in the previous section.  Again, $I_{\omega,0}$ and the terms involving $F_0(\alpha k)$ and $F_0(\alpha^2 k)$ are dropped when $s > 0$.   A general expansion follows
\begin{align*}
q(x,t)  &= \frac{1}{2 \pi} \int_{C} e^{iks} F_0(k) \frac{dk}{ik} - \frac{\alpha^{-1}}{2 \pi} \int_{\partial D^+_2} e^{iks} F_0( \alpha k) \frac{dk}{ik}  - \frac{\alpha^{-2}}{2 \pi} \int_{\partial D^+_1} e^{iks} F_0( \alpha^2 k) \frac{dk}{ik}  \\
&+ \sum_{x_i} [q_o(x_i)] I_{\omega,0}(x-x_i,t) + q_o(0) \left(I_{\omega,0}(x,t) + \alpha^{-1} I_{\omega,0,2}(x,t) + \alpha^{-2} I_{\omega,0,1}(x,t) \right)\\
& +   \bigo(|x-s|^{1/2} + t^{1/6}).
\end{align*}
Here the integrals on $\partial D^+_1$ and $\partial D^+_2$ should be dropped when $s > 0$.

For reasons that are made clear below, we require another iteration of integration by parts for $s = 0$.  In the case that the first derivative of $q_o$ has no discontinuities we have
\begin{align*}
\hat q_o(k) = \frac{q_o(0)}{ik} + \frac{q'_o(0)}{(ik)^2} + \frac{F_1(k)}{(ik)^2}, \quad F_1(k) = \int_{0}^{\infty} e^{-iks} q_o''(s) ds.
\end{align*}
Then
\begin{align*}
q(x,t) &= q_o(0) \left(I_{\omega,0}(x,t) - \alpha^{-1} I_{\omega,0,2}(x,t) - \alpha^{-2} I_{\omega,0,1}(x,t) \right) \\
&+ q_o'(0) \left(I_{\omega,1}(x,t) - \alpha^{-2} I_{\omega,1,2}(x,t) - \alpha^{-4} I_{\omega,1,1}(x,t) \right)\\
&+ \frac{1}{2 \pi} \int_{C} e^{ikx-i\omega(k)t} \frac{F_1(k)}{(ik)^2} dk - \frac{\alpha^{-2}}{2 \pi} \int_{\partial D^+_2} e^{ikx-i\omega(k)t} \frac{F_1(\alpha k)}{(ik)^2} dk - \frac{\alpha^{-4}}{2 \pi} \int_{\partial D^+_1} e^{ikx-i\omega(k)t} \frac{F_1(\alpha^2 k)}{(ik)^2} dk,
\end{align*}
and
\begin{align*}
q(x,t) &=  q_o(0) \left(I_{\omega,0}(x,t) - \alpha^{-1} I_{\omega,0,2}(x,t) - \alpha^{-2} I_{\omega,0,1}(x,t) \right) \\
&+ q_o'(0) \left(I_{\omega,1}(x,t) - \alpha^{-2} I_{\omega,1,2}(x,t) - \alpha^{-4} I_{\omega,1,1}(x,t) \right)\\
&+ \frac{1}{2 \pi} \int_{C}(1+ikx) \frac{F_1(k)}{(ik)^2} dk - \frac{\alpha^{-2}}{2 \pi} \int_{\partial D^+_2} (1+ikx)\frac{F_1(\alpha k)}{(ik)^2} dk \\
&- \frac{\alpha^{-4}}{2 \pi} \int_{\partial D^+_1}(1+ikx)\frac{F_1(\alpha^2 k)}{(ik)^2} dk + \bigo\left( x^{3/2} + t^{1/2} \right).
\end{align*}
From this it should be clear how to treat the case of multiple discontinuities in $q_o$ and $q_o'$.

%%% Local Variables:
%%% mode: latex
%%% TeX-master: "gibbs-ibvp"
%%% End:

\section{IBVP with zero initial data}
\label{sec:zeroIC}

In this section we treat the case where the initial data for the IBVP vanishes identically. Linearity allows us to combine the results from this section with that of the previous section to produce a full characterization of the solution near the boundary under Assumption~\ref{Assume:Data-0}.  Furthermore, following ideas from Appendix~\ref{a:SolutionForm} it suffices to treat the case where the boundary data is in $H^1([0,T])$:  Any other discontinuities can be added through linearity.  For zero initial data there are three relevant components of the analysis of this solution formula:
\begin{enumerate}
\item the behavior of $q$ near $x = 0$ for $t > 0$,
\item the behavior of $q$ near $(x,t) = (0,0)$,
\item the behavior of $q$ near $(x,t) = (s,0)$  for $0 < s < T$.
\end{enumerate}

\subsection{Linear Schr\"odinger}

With zero initial data the solution of \eref{e:LS} is simply given by ($\omega(k) = k^2$)
\begin{align}\label{e:ls-zeroIC}
q(x,t) = \frac{1}{2\pi} \int_{\partial D_1^+} e^{ikx-i\omega(k)t} 2k \tilde g_0(-\omega(k),t) dk.
\end{align}
We integrate $\tilde g_0(k,t)$ by parts.  This gives
\begin{align*}
\tilde g_0(-\omega(k),t) &= \frac{g_0(t) e^{i\omega(k)t}-g_0(0)}{i\omega(k)} - \frac{G_{0,0}(k)}{i\omega(k)},\\
G_{0,0}(k) &=  \int_{0}^t e^{i \omega(k) s} g_0'(s) ds.
\end{align*}
Then
\begin{align*}
q(x,t) = - 2 g_0(0) I_{\omega,0,1}(x,t) - \frac{1}{\pi}\int_{\partial D_1^+} e^{ikx-i\omega(k)t} G_{0,0}(k) \frac{dk}{ik},
\end{align*}
because the term involving $g_0(t)$ vanishes by Jordan's Lemma. Furthermore, all of these functions are continuous up to $x = 0$.  When considering \eref{e:ibp-LS} we see that the contribution from $I_{\omega,0,1}$ will cancel if these two solutions are added and the first compatibility condition holds: $q_o(0) = g_0(0)$.  We then appeal to Lemmas~\ref{Lemma:Expansion} and \ref{Lemma:Sufficient} to derive the expansion near $(s,\tau)$,
\begin{align}\label{e:ls-zeroIC-ex}
q(x,t) = - 2 g_0(0) I_{\omega,0,1}(x,t) - \frac{1}{\pi}\int_{\partial D_1^+} e^{iks-i\omega(k)\tau} G_{0,0}(k) \frac{dk}{ik} + \bigo( |x-s|^{1/2} + |t-\tau|^{1/4}).
\end{align}
This is the correct form for the solution when $s =0$, $\tau = 0$.  This formula is now further examined in the remaining regimes discussed above.  For $s > 0$ and $\tau = 0$, $q(x,t) = \bigo(|t|^{1/4})$.  For $s = 0$, $\tau > 0$ we claim
\begin{align*}
q(x,t) &= - 2 g_0(0) I_{\omega,0,1}(x,t) - \frac{1}{\pi}\int_{\partial D_1^+} e^{-i\omega(k)\tau} G_{0,0}(k) \frac{dk}{ik} + \bigo( |x|^{1/2} + |t-\tau|^{1/4}) \\
&= g_0(\tau) +  \bigo( |x|^{1/2} + |t-\tau|^{1/4}).
\end{align*}
Indeed,
\begin{align*}
- 2 g_0(0) I_{\omega,0,1}(x,t) &- \frac{1}{\pi}\int_{\partial D_1^+} e^{-i\omega(k)\tau} G_{0,0}(k) \frac{dk}{ik} \\
&= ( 2 g_0(0) I_{\omega,0,1}(0,\tau) - 2 g_0(0) I_{\omega,0,1}(x,t)) - 2 g_0(0) I_{\omega,0,1}(0,\tau) \\
&- \frac{1}{\pi}\int_{\partial D_1^+} e^{-i\omega(k)\tau} G_{0,0}(k) \frac{dk}{ik}\\
&= 2 g_0(0) (I_{\omega,0,1}(0,\tau) - I_{\omega,0,1}(x,t)) + g_0(\tau).
\end{align*}
This follows from:
\begin{lemma}\label{l:g}
  For $0 < \tau < T$ and $g_0 \in H^1([0,T])$,
  \begin{align*}
  g_0(\tau) = - 2 g_0(0) I_{\omega,0,1}(0,\tau) - \frac{1}{\pi}\int_{\partial D_1^+} e^{-i\omega(k)\tau} G_{0,0}(k) \frac{dk}{ik}.
  \end{align*}
\end{lemma}
\begin{proof}
  First,  it follows that $I_{\omega,0,1}(0,\tau) = -1/2$ for $\tau > 0$. Then it suffices to show
  \begin{align*}
\int_0^\tau g'_0(s) ds = -\frac{1}{\pi}\int_{\partial D_1^+} e^{-i\omega(k)s} G_{0,0}(k) \frac{dk}{ik}.
  \end{align*}
   Using $l = k^2 = \omega(k)$, for a.e. $s \in [0,\tau]$
\begin{align*}
g'(s) &= \frac{1}{2 \pi} \int_{\mathbb R} e^{-isl} \int_0^\tau e^{is'l} g'(s') ds' dl\\
&= \frac{1}{2 \pi} \int_{\partial D^+_1} 2k e^{-i\omega(k)s} \int_0^\tau e^{i\omega(k)s'} g'(s') ds' dk\\
&= \frac{1}{2\pi} \int_{\partial D_1^+}  e^{-i\omega(k)s} 2k G_{0,0}(k) dk.
\end{align*}
We need to justify integrating this expression with respect to $s$ and interchanging the order of integration.  Let $\Gamma_R = B(0,R) \cap \partial D_1^+$ and we have
\begin{align*}
\int_0^\tau g'(s) ds &= \int_0^\tau \lim_{R \rightarrow \infty} \frac{1}{2\pi} \int_{\Gamma_R}  e^{-i\omega(k)s} 2k G_{0,0}(k) dk ds\\
&= \lim_{R \rightarrow \infty}\int_0^\tau  \frac{1}{2\pi} \int_{\Gamma_R}  e^{-i\omega(k)s} 2k G_{0,0}(k) dk ds
\end{align*}
by the dominated convergence theorem.   Now, because we have finite domains for the integration of bounded functions we can interchange:
\begin{align*}
\int_0^\tau g'(s) ds  &= \lim_{R \rightarrow \infty} \int_{\Gamma_R} \int_0^\tau  \frac{1}{2\pi}   e^{-i\omega(k)s} 2k G_{0,0}(k) dk ds\\
&= \lim_{R \rightarrow \infty}\frac{1}{2\pi} \int_{\Gamma_R} [e^{-i\omega(k)\tau}-1] \frac{2k}{-i \omega(k)} G_{0,0}(k) dk\\
&= -\lim_{R \rightarrow \infty} \frac{1}{\pi}\int_{\Gamma_R} e^{-i\omega(k)\tau} \frac{1}{ik} G_{0,0}(k) dk \\
&+\lim_{R \rightarrow \infty} \frac{1}{\pi}\int_{\Gamma_R} \frac{1}{ik} G_{0,0}(k) dk \\
&= - \frac{1}{\pi}\int_{\partial D_1^+} e^{-i\omega(k)\tau} \frac{1}{ik} G_{0,0}(k) dk,
\end{align*}
because the integral in the second-to-last  line vanishes from Jordan's Lemma.
\end{proof}

Then \eqref{e:ls-zeroIC-ex} follows because $I_{\omega,0,1}(x,t)$ is a smooth function of $(x,t)$ for $t > 0$. So, \eqref{e:ls-zeroIC-ex} is the expansion about $(s,\tau)$ for any choice of $(s,\tau)$ in $\overline{\mathbb R^+ \times (0,T)}$, including $(s,\tau) = (0,0)$.  As the calculations get more involved in the following sections, we skip calculations along the lines of Lemma~\ref{l:g}.

\subsection{Airy 1}

In the case of \eref{e:Airy1} with zero initial data we have ($\omega(k) = -k^3$)
\begin{align*}
q(x,t) &= -\frac{1}{2\pi}\int_{\partial D^+}3k^2e^{ikx-i\omega(k)t}\tilde{g}_0(-\omega(k),t) dk.
\end{align*}
Integration by parts gives the expansion
\begin{align}\label{e:A1-zeroIC-ex}
q(x,t) = - 3 g_0(0) I_{\omega,0,1}(x,t) - \frac{1}{2\pi}\int_{\partial D_1^+} e^{iks-i\omega(k)\tau} G_{0,0}(k) \frac{dk}{ik} + \bigo(|x-s|^{1/2} + |t-\tau|^{1/6}).
\end{align}
This right-hand side is easily seen to be $\bigo(|x-s|^{1/2} + |t-\tau|^{1/6})$ when $s > 0$ and $\tau = 0$.  Additionally, for $s = 0$ and $\tau > 0$ it follows in a similar manner to Lemma~\ref{l:g} that 
\begin{align*}
q(x,t) &= g_0(\tau) +  \bigo( |x|^{1/2} + |t-\tau|^{1/6}).
\end{align*}
As in the previous case \eqref{e:A1-zeroIC-ex} is the appropriate expansion about $(s,\tau)$ for any choice of $(s,\tau)$ in $\overline{\mathbb R^+ \times \mathbb (0,T)}$, including $(s,\tau) = (0,0)$.

\subsection{Airy 2}

We consider the more interesting case of \eqref{e:Airy2}.  Here $\omega(k) = k^3$ and the solution is given by
\begin{align*}
q(x,t) =  - \frac{1}{2\pi}\int_{\partial D_2^+} e^{ikx-i\omega(k)t}\left( (\alpha^2-1)k^2\tilde{g}_0(-\omega(k),t) - i(\alpha-1) k\tilde{g}_1(-\omega(k),t)  \right)dk \\
- \frac{1}{2\pi}\int_{\partial D_1^+} e^{ikx-i\omega(k)t} \left ((\alpha-1) k^2\tilde{g}_0(-\omega(k),t) - i(\alpha^2-1) k\tilde{g}_1(-\omega(k),t) \right) dk.
\end{align*}
We integrate both $\tilde g_0$ and $\tilde g_1$ by parts
\begin{align*}
\tilde g_0(-\omega(k),t) &= \frac{g_0(t) e^{i\omega(k)t}-g_0(0)}{i\omega(k)} -  \frac{g_0'(t) e^{i\omega(k)t}-g_0'(0)}{(i\omega(k))^2} + \frac{G_{0,1}(k)}{(i\omega(k))^2},\\
\tilde g_1(-\omega(k),t) &= \frac{g_1(t) e^{i\omega(k)t}-g_1(0)}{i\omega(k)} - \frac{G_{1,0}(k)}{i\omega(k)},\\
G_{1,0}(k) &=  \int_{0}^t e^{i \omega(k) s} g_1'(s) ds, \quad G_{0,1}(k) =  \int_{0}^t e^{i \omega(k) s} g_0''(s) ds.
\end{align*}
We then see that
\begin{align}\label{e:I1}
\begin{split}
\mathcal I_1(x,t) &:= \frac{1}{2\pi}\int_{\partial D_2^+} e^{ikx-i\omega(k)t} k^2 \tilde{g}_0(-\omega(k),t) dk \\
&= \frac{1}{2\pi}\int_{\partial D_2^+} e^{ikx-i\omega(k)t} \left(\frac{g_0(t) e^{i\omega(k)t}-g_0(0)}{ik} -  \frac{g_0'(t) e^{i\omega(k)t}-g_0'(0)}{(ik)^2k^2} + \frac{G_{0,1}(k)}{(ik)^2k^2}\right) dk.
\end{split}
\end{align}
Terms with the factor $e^{i\omega(k)t}$ vanish by Jordan's lemma so that
\begin{align*}
\mathcal I_1(x,t) &= -\frac{1}{2\pi}\int_{\partial D_2^+} e^{ikx-i\omega(k)t} \left(\frac{g_0(0)}{ik} -  \frac{g_0'(0)}{(ik)^2k^2} - \frac{G_{0,1}(k)}{(ik)^2k^2}\right) dk\\
&= -g_0(0) I_{\omega,0,2}(x,t)  \\
&- \frac{1}{2\pi}\int_{\partial D_2^+} (1 + ik(x-s) +  k^2 (x-s)^2 - i k^3 (x-s)^3 - i\omega(k) (t-\tau)) e^{iks-i\omega(k)\tau}\frac{g_0'(0) -G_{0,1}(k)}{(ik)^2k^2} dk \\
&+ \bigo( |x-s|^{7/2} + |t-\tau|^{7/6})\\
&= -g_0(0) I_{\omega,0,2}(x,t)  + \frac{1}{2\pi}\int_{\partial D_2^+}(1 + ik(x-s))\frac{g_0'(0) -G_{0,1}(k)}{(ik)^2k^2} dk + \bigo(|x-s|^{2} + |t-\tau|).
\end{align*}
We only need to keep the terms involving $(x-s)$. Next, we consider
\begin{align}\label{e:I2}
\begin{split}
\mathcal I_2(x,t) &:= \frac{1}{2\pi}\int_{\partial D_2^+} e^{ikx-i\omega(k)t} k \tilde{g}_1(-\omega(k),t) dk \\
&= \frac{i}{2\pi}\int_{\partial D_2^+} e^{ikx-i\omega(k)t} \left(\frac{g_1(t) e^{i\omega(k)t}-g_1(0)}{(ik)^2} -  \frac{G_{1,0}(k)}{(ik)^2}\right) dk\\
& =  i g_1(0) I_{\omega,1,2}(x,t) - \frac{i}{2\pi}  \int_{\partial D_2^+}(1+ik(x-s)) e^{iks-i\omega(k)\tau}\frac{G_{1,0}(k)}{(ik)^2} dk \\
&+ \bigo(|x-s|^{3/2} + |t-\tau|^{1/2}).
\end{split}
\end{align}
Combining all of this with the integrals on $\partial D_1^+$ we find
\begin{align}\label{e:A2-zeroIC-ex}
\begin{split}
q(x,t) &= g_0(0) \left((\alpha^2-1)I_{\omega,0,2}(x,t)+(\alpha-1) I_{\omega,0,1}(x,t) \right)\\
  &-g_1(0) \left( (1-\alpha)I_{\omega,1,2}(x,t) + (1-\alpha^2)I_{\omega,1,1}(x,t) \right)\\
&- \frac{1- \alpha^2}{2 \pi}  \int_{\partial D_2^+}(1+ik(x-s))e^{iks-i\omega(k)\tau}\frac{g_0'(0) -G_{0,1}(k)}{(ik)^2k^2} dk \\
&+ \frac{\alpha -1}{2 \pi}  \int_{\partial D_2^+} (1+ik(x-s))e^{iks-i\omega(k)\tau}\frac{G_{1,0}(k)}{(ik)^2} dk\\
&  - \frac{1- \alpha}{2 \pi} \int_{\partial D_1^+}(1+ik(x-s))e^{iks-i\omega(k)\tau}\frac{g_0'(0) -G_{0,1}(k)}{(ik)^2k^2} dk  \\
&+ \frac{\alpha^2 -1}{2\pi} \int_{\partial D_2^+} (1+ik(x-s))e^{iks-i\omega(k)\tau}\frac{G_{1,0}(k)}{(ik)^2} dk\\
& + \bigo(|x-s|^{3/2} + |t-\tau|^{1/2}).
\end{split}
\end{align}
If $s > 0$ and $\tau = 0$ then all integrals along $\partial D_i^+$ for $i =1,2$ vanish identically and $q(x,t) = \bigo (|x-s|^{3/2} + |t|^{1/2})$.  To analyze the expression when  $s = 0$ and $\tau > 0$, we consider
\begin{align*}
\mathcal L_0(\tau) &:= g_0(0) \left((\alpha^2-1)I_{\omega,0,2}(0,\tau)+(\alpha-1) I_{\omega,0,1}(0,\tau) \right)\\
&+ g_1(0) \left( (1-\alpha)I_{\omega,1,2}(0,\tau) + (1-\alpha^2)I_{\omega,1,1}(0,\tau) \right)\\
&- \frac{1- \alpha^2}{2 \pi}  \int_{\partial D_2^+}e^{-i\omega(k)\tau}\frac{g_0'(0) -G_{0,1}(k)}{(ik)^2k^2} dk + \frac{\alpha -1}{2 \pi}  \int_{\partial D_2^+} e^{-i\omega(k)\tau}\frac{G_{1,0}(k)}{(ik)^2} dk\\
&  - \frac{1- \alpha}{2 \pi} \int_{\partial D_1^+}e^{-i\omega(k)\tau}\frac{g_0'(0) -G_{0,1}(k)}{(ik)^2k^2} dk  + \frac{\alpha^2 -1}{2\pi} \int_{\partial D_1^+} e^{-i\omega(k)\tau}\frac{G_{1,0}(k)}{(ik)^2} dk.\\
\end{align*}
Because multiplication by $\alpha^{-1}$ takes $\partial D_2^+$ to $\partial D_1^+$, and $G_{i,j}(\alpha k) = G_{i,j}(k)$ we find
\begin{align*}
\frac{1- \alpha}{2 \pi} \int_{\partial D_1^+}e^{-i\omega(k)\tau}\frac{g_0'(0) -G_{0,1}(k)}{(ik)^2k^2} dk = \frac{1- \alpha}{2 \pi} \int_{\partial D_2^+}e^{-i\omega(k)\tau}\frac{g_0'(0) -G_{0,1}(k)}{(ik)^2k^2} dk, \\
\frac{\alpha^2 -1}{2\pi} \int_{\partial D_1^+} e^{-i\omega(k)\tau}\frac{G_{1,0}(k)}{(ik)^2} dk = \frac{1 - \alpha}{2\pi} \int_{\partial D_2^+} e^{-i\omega(k)\tau}\frac{G_{1,0}(k)}{(ik)^2} dk.
\end{align*}
Thus, the terms involving $G_{1,0}(k)$ and $I_{\omega,1,j}$ vanish identically and it can be shown that $\mathcal L_0(\tau) = g_0(\tau)$.  Then we consider a term that resembles differentiation in $x$
\begin{align*}
\mathcal L_1(\tau) &:= g_0(0) \left((\alpha^2-1)I_{\omega,-1,2}(0,\tau)+(\alpha-1) I_{\omega,-1,1}(0,\tau) \right)\\
&+ g_1(0) \left( (1-\alpha)I_{\omega,1,2}(0,\tau) + (1-\alpha^2)I_{\omega,1,1}(0,\tau) \right)\\
&- \frac{1- \alpha^2}{2 \pi}  \int_{\partial D_2^+}e^{-i\omega(k)\tau}\frac{g_0'(0) -G_{0,1}(k)}{(ik)k^2} dk + \frac{\alpha -1}{2 \pi}  \int_{\partial D_2^+} e^{-i\omega(k)\tau}\frac{G_{1,0}(k)}{ik} dk\\
&  - \frac{1- \alpha}{2 \pi} \int_{\partial D_1^+}e^{-i\omega(k)\tau}\frac{g_0'(0) -G_{0,1}(k)}{(ik)k^2} dk  + \frac{\alpha^2 -1}{2\pi} \int_{\partial D_1^+} e^{-i\omega(k)\tau}\frac{G_{1,0}(k)}{ik} dk.
\end{align*}
We use
\begin{align}
\frac{1- \alpha}{2 \pi} \int_{\partial D_1^+}e^{-i\omega(k)\tau}\frac{g_0'(0) -G_{0,1}(k)}{(ik)k^2} dk = \frac{\alpha^2- 1}{2 \pi} \int_{\partial D_2^+}e^{-i\omega(k)\tau}\frac{g_0'(0) -G_{0,1}(k)}{(ik)k^2} dk,\notag\\
\frac{\alpha^2 -1}{2\pi} \int_{\partial D_1^+} e^{-i\omega(k)\tau}\frac{G_{1,0}(k)}{ik} dk = \frac{\alpha^2 -1}{2\pi} \int_{\partial D_2^+} e^{-i\omega(k)\tau}\frac{G_{1,0}(k)}{ik} dk,\label{D1-D2}
\end{align}
to see that all terms involving $g_0$ cancel identically. It then can be shown that
\begin{align*}
\mathcal L_1(\tau) = g_1(\tau),
\end{align*}
and finally
\begin{align*}
q(x,t) = g_1(\tau) + x g_1(\tau) + \bigo(|x|^{3/2} + |t-\tau|^{1/2}),
\end{align*}
as expected.  Again, \eqref{e:A2-zeroIC-ex} is the appropriate expansion about $(s,\tau)$ for any choice of $(s,\tau)$ in $\overline{\mathbb R^+ \times (0,T)}$, including $(s,\tau) = (0,0)$.

%%% Local Variables:
%%% mode: latex
%%% TeX-master: "gibbs-ibvp.tex"
%%% End:

%\input s5higherorder
\section{Higher-order theory and decay of the spectral data}
\label{sec:spectraldata}

If the initial and boundary data are compatible in the sense that $q_o(0) = g_o(x)$ it is straightforward to check in the examples considered that the terms involving $I_{\omega,0,j}(x,t)$ drop out of the solution formula after integration by parts. The expressions from Section~\ref{sec:zeroIC} are added to those from Section~\ref{sec:zeroBC} to see this.  Furthermore, in the case of \eref{e:Airy2} if $q_o'(0) = g_1(0)$ then the terms $I_{\omega,1,j}$ drop out.  This is related to the fact that smoothness of the data plus higher-order compatibility at the corner $(x,t) = (0,0)$ forces the integrands in \eqref{e:form} to decay more rapidly.  Specifically, it is clear that the expressions for $\mathcal I_1$ and $\mathcal I_2$ (see \eref{e:I1} and \eref{e:I2}) once $I_{\omega,m,j}$ are removed have integrands that decay faster.   Understanding this behavior is important for many reasons, one of which is numerical evaluation.

We trust that our example here is enough to demonstrate the relevant behavior when the initial and boundary data are compatible.  We focus on \eref{e:Airy2} and apply repeated integration by parts.  We only write the terms that involve the functions $I_{\omega,m,j}$.  It is clear by using $I_{\omega,m}(x,t) = I_{\omega,m,1}(x,t) + I_{\omega,m,2}(x,t)$ that
\begin{align}\label{e:qo}
q|_{g_j \equiv 0}(x,t) = \sum_{i=0}^\ell q^{(i)}(0) \left( (1 - \alpha^{-1-i}) I_{\omega,i,1}(x,t) + (1 - \alpha^{-2-2i}) I_{\omega,i,2}(x,t) \right) + E_{g_j \equiv 0}(x,t).
\end{align}
Here $E_{g_j\equiv 0}$ represents components of the solution not expressed in terms of $I_{\omega,m,j}$.  Next using that $\alpha^2 = \alpha^{-1}$ and $\alpha = \alpha^{-2}$
\begin{align}\label{e:g0}
q|_{q_o \equiv 0}(x,t) &= \sum_{j=0}^{\ell} g_0^{(j)}(0) (  (\alpha^{-1}-1) I_{\omega,3j,1}(x,t) + (\alpha^{-2}-1) I_{\omega,3j,2}(x,t))\\
&+ \sum_{j=0}^{\ell} g_1^{(j)}(0) ( (\alpha^{-2}-1) I_{\omega,3j+1,1}(x,t) + (\alpha^{-1}-1) I_{\omega,3j+1,2}(x,t)) + E|_{q_o\equiv 0}(x,t). \label{e:g1}
\end{align}
We consider cancellations in the sum $q|_{q_o\equiv 0} + q|_{g_j\equiv 0}$.  Now, if $i = 3j$ then
\begin{align*}
(1 - \alpha^{-1-i}) I_{\omega,i,1}(x,t) + (1 - \alpha^{-2-2i}) I_{\omega,i,2}(x,t) = (1 - \alpha^{-1}) I_{\omega,3j,1}(x,t) + (1 - \alpha^{-2}) I_{\omega,3j,2}(x,t).
\end{align*}
 If $q_o^{3j}(0) = g_0^{(j)}(0)$ one term in the sums in \eref{e:qo} and \eref{e:g0} cancel.  Now, if $i = 3j+1$ a similar cancellation occurs if $q_o^{3j+1}(0) = g_1^{(j)}(0)$.  Thus, it remains to consider $i = 3j +2$.  In this case, a simple calculation reveals $\alpha^{-1 - (3j+2)} = \alpha^{-2 -2(3j+2)} = 1$ and cancellation of this term requires no additional conditions on the initial/boundary data.  What we have displayed is the following.

\begin{proposition}
Assume $q_o \in H^{m}(\mathbb R^+)$  and $g_j \in H^{\lceil (m-j)/n \rceil}(\mathbb R)$ for $j = 0, \ldots,N(n)-1$.  Further, assume the compatibility conditions hold up to order $m$.  Then the spectral data, \emph{i.e} the integrand $\mathcal F$ of \eref{e:form} at $x = t = 0$, can be written so that it satisfies
\begin{align*}
\mathcal F(\cdot) (1 + |\cdot|)^m \in L^2(\partial D).
\end{align*}
\end{proposition}

We do not present the details here but to obtain an asymptotic expansion for $q(x,t)$ when discontinuities exist in higher-order derivatives, one applies Lemma~\ref{Lemma:Expansion} (after the cancellation of appropriate terms involving $I_{\omega,i,j}$) to expand terms of the form
\begin{align*}
\int_{\partial D_i^+} \frac{F_j(k)}{(ik)^{j+1}} dk, \quad \int_{\partial D_i^+} \frac{G_{j,\ell}(k)}{(i\omega(k))^j k^m} dk,
\end{align*}
which result from integration by parts.

%%% Local Variables:
%%% mode: latex
%%% TeX-master: "gibbs-ibvp.tex"
%%% End:

\section{Example solutions of IBVPs with general corner singularities}
\label{sec:examples}

We now combine the results of the previous sections and we discuss the behavior of the solutions
of the IBVP when the ICs and BCs are both non-zero, but one of the compatibility conditions is violated.  We note that because of the expansions above, the dominant behavior of the solution near any discontinuity in the data is given in terms of the special functions $I_{\omega,m,j}(x,t)$ and we focus on plotting this dominant behavior.

A few words should be said about computing $I_{\omega,m,j}(x,t)$.  When using the steepest method for integrals as in Theorem~\ref{Thm:Kernel} (again see \cite{GinoTomIVP} for details) the path of steepest descent can be approximated and a numerical quadrature routine applied on this approximate contour.  With some care to \br scale contours \er appropriately near the stationary phase point as, the method is provably accurate for all values of the parameters.  We refer the reader to a discussion of this in \cite{TrogdonUnifiedNumerics} and in \cite{GinoTomIVP}.  In what follows, we use Clenshaw--Curtis quadrature \cite{clenshaw} on piecewise affine contours which is implemented in {\tt RHPackage} \cite{RHPackage} and we are able to approximate any one of the functions $I_{\omega,m,j}(x,t)$ well, even as $x \to \infty$ or $t \downarrow 0$.

%%%%%%%%%%%%%%%%%%%%%%%%%%%%%%%%%%%%%%%%%%%%%%%%%%%%%%%%%%%%%%%%%%%%%%%%%%%%%%%%%%%%%%%%%%%%%
\subsection{Linear Schr\"odinger}
\newcommand{\qloc}{q_{\text{loc}}}

\begin{figure}[htp]
\subfigure[]{\includegraphics[width=.49\linewidth]{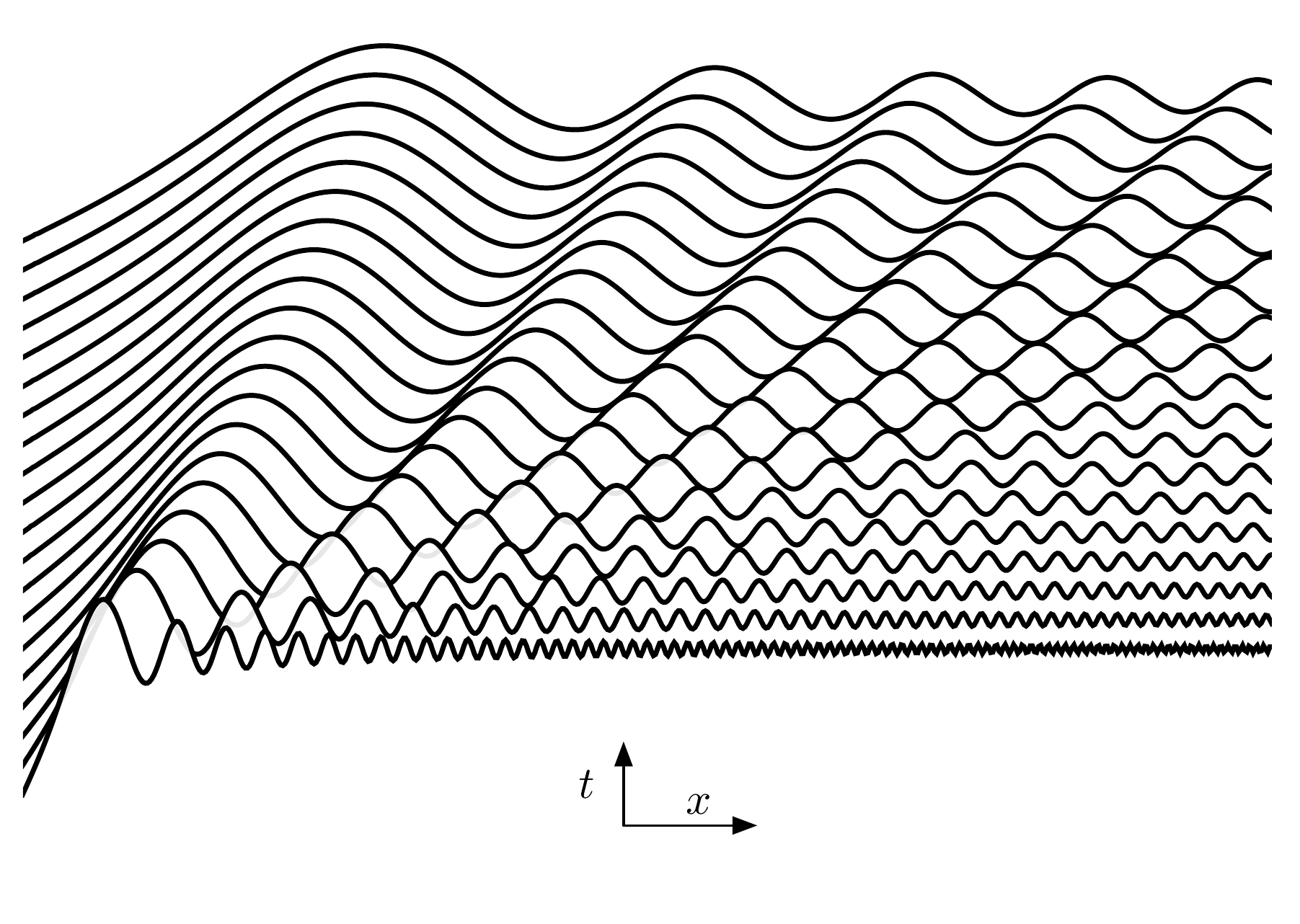}}
\subfigure[]{\includegraphics[width=.49\linewidth]{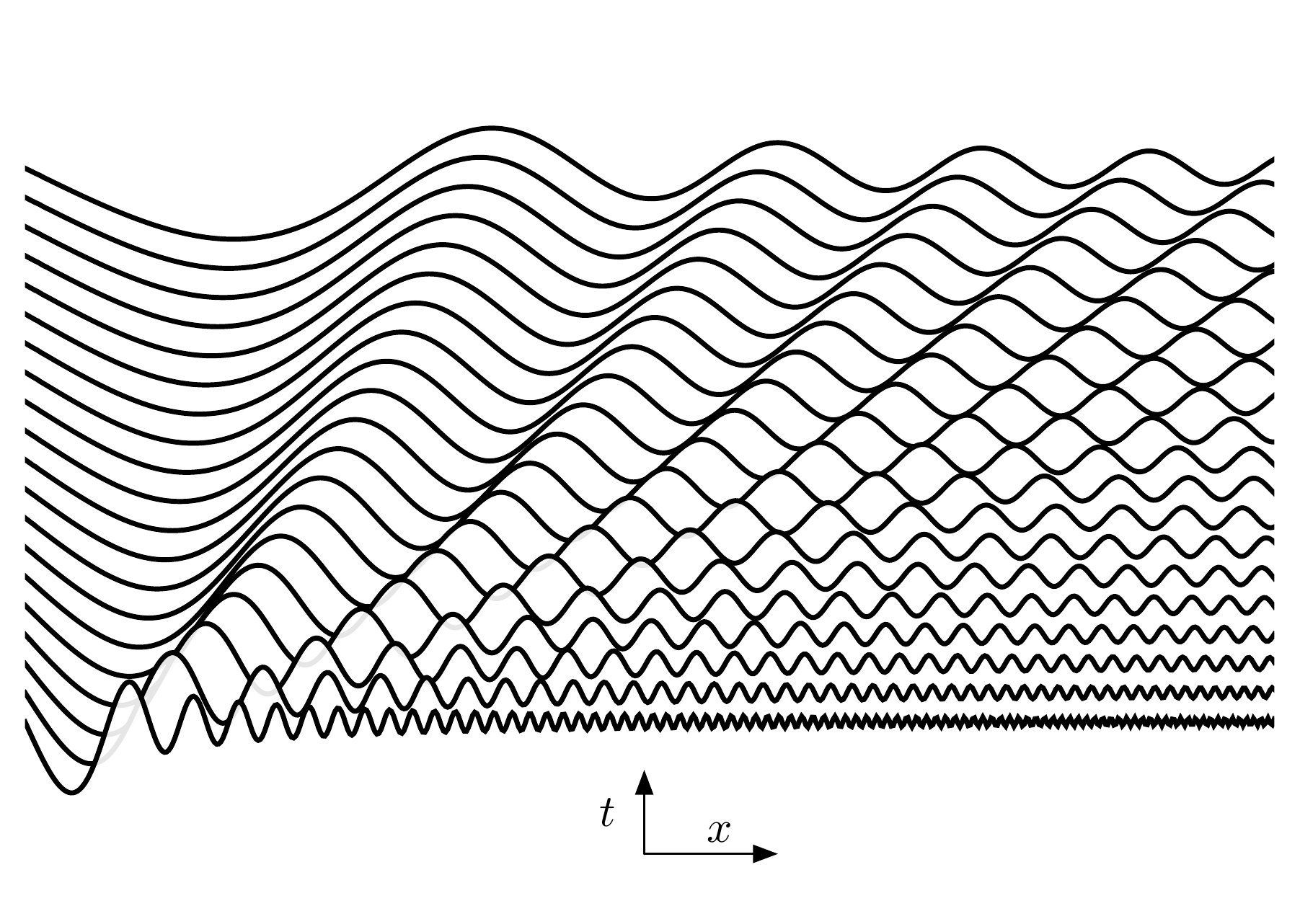}}
\subfigure[]{\includegraphics[width=.49\linewidth]{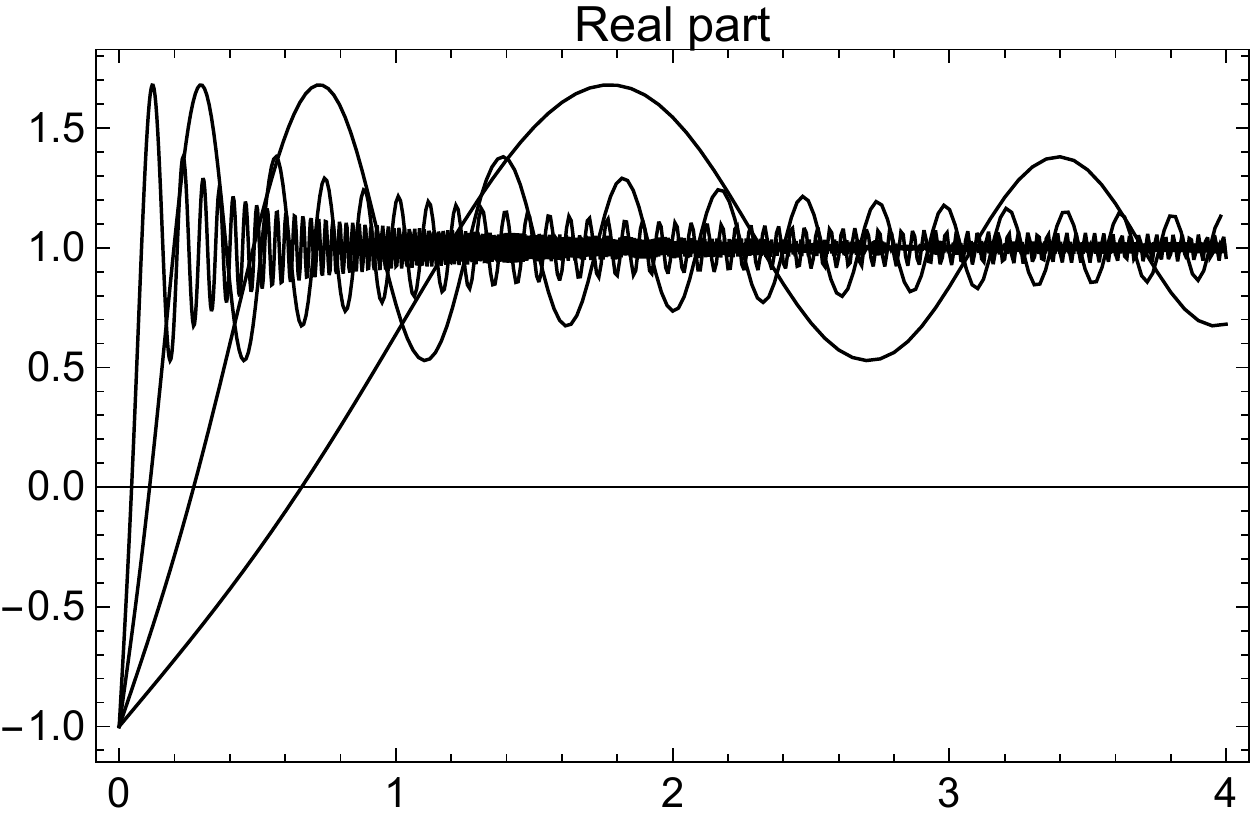}}
\subfigure[]{\includegraphics[width=.49\linewidth]{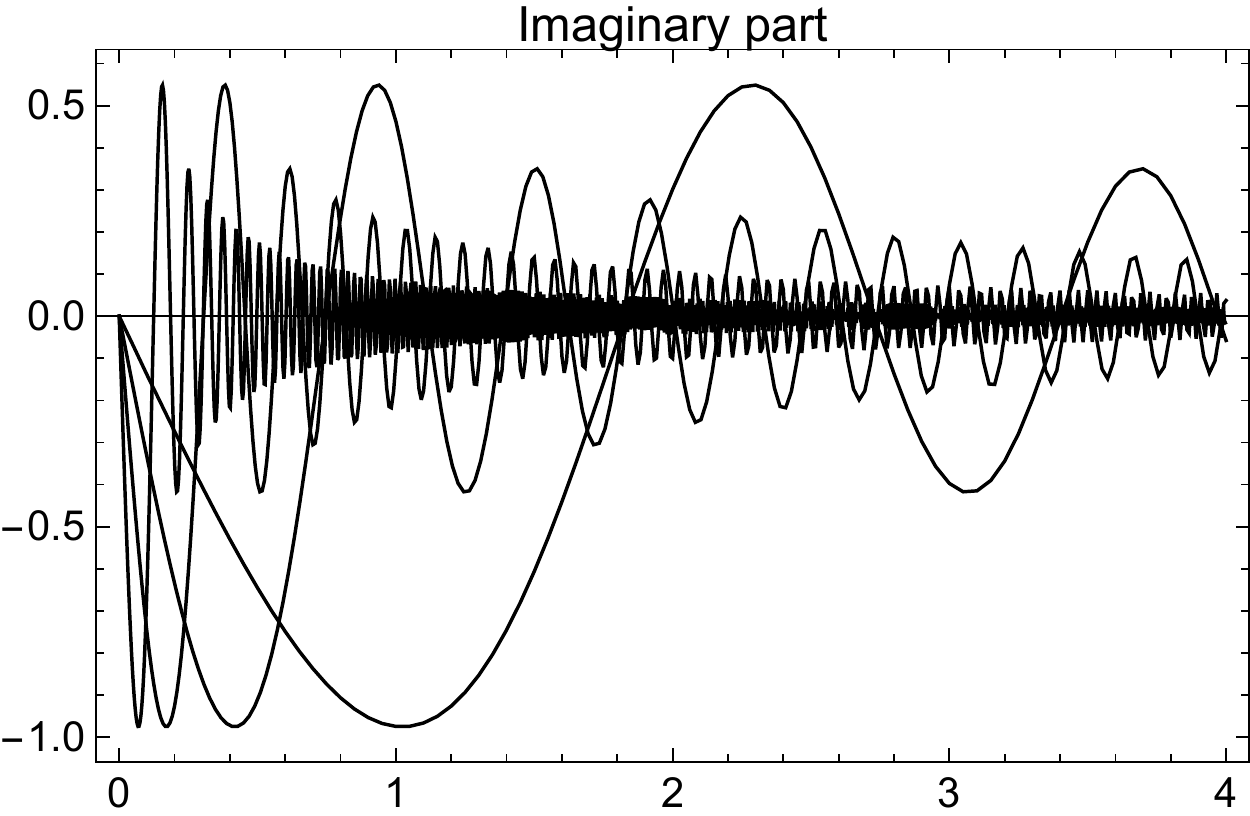}}
\caption{\label{f:LS-loc} Plots of $\qloc(x,t)$ for the linear Schr\"odinger equation in the concrete case $q_o(0) = 1$ and $g_0(0) = -1$. (a) The time evolution of $\re \qloc(x,t)$ up to $t = 2$.  (b) The time evolution of $\im \qloc(x,t)$ up to $t = 2$.  (c) An examination of $\re \qloc(x,t)$ as $t \downarrow 0$ for $t = 1/20(1/6)^j$, $j = 0,1,2,3$. It is clear that the solution is limiting to $\qloc(x,t) = 1$ for $x > 0$ and satisfies $\qloc(0,t) = -1$ for all $t$.  (d) An examination of $\im \qloc(x,t)$ as $t \downarrow 0$ for $t = 1/20(1/6)^j$, $j = 0,1,2,3$.}
\end{figure}

If we were to examine the solution of \eqref{e:LS} near a corner singularity with $\omega(k) = k^2$ we would be led to the expansion
\begin{align*}
q(x,t) = 2 (q_o(0) - g_0(0)) I_{\omega,0,1}(x,t) + C+ \bigo(|x|^{1/2} + |t|^{1/4}).
\end{align*}
The constant $C$ is given in terms of integrals of $F_0$ and $G_{0,0}$ but it can be found by other reasoning.  For example, if we set $x = 0$ and let $t \downarrow 0$ then $\lim_{t \downarrow 0} q(x,t) = g_0(0)$.  It follows from Theorem~\ref{Thm:Kernel} that  $\lim_{t \downarrow 0}I_{\omega,0,1}(x,t) = 0$ for $x> 0$ so that $C = q_o(0)$ and the solution is
\begin{align*}
q(x,t) &= q_{\text{loc}}(x,t) + \bigo(|x|^{1/2} + |t|^{1/4}),\\
\qloc(x,t) &= -2 g_0(0) I_{\omega,0,1}(x,t) + 2q_o(0) \left(I_{\omega,0,1}(x,t)+\frac{1}{2}\right).
\end{align*}
A concrete case is $q_o(0) = 1$ and $g_0(0) = -1$ and we explore $\qloc(x,t)$ in Figure~\ref{f:LS-loc}.

%%%%%%%%%%%%%%%%%%%%%%%%%%%%%%%%%%%%%%%%%%%%%%%%%%%%%%%%%%%%%%%%%%%%%%%%%%%%%%%%%%%%%%%%%%%%%
\subsection{Airy 1}

\begin{figure}[htp]
\subfigure[]{\includegraphics[width=.49\linewidth]{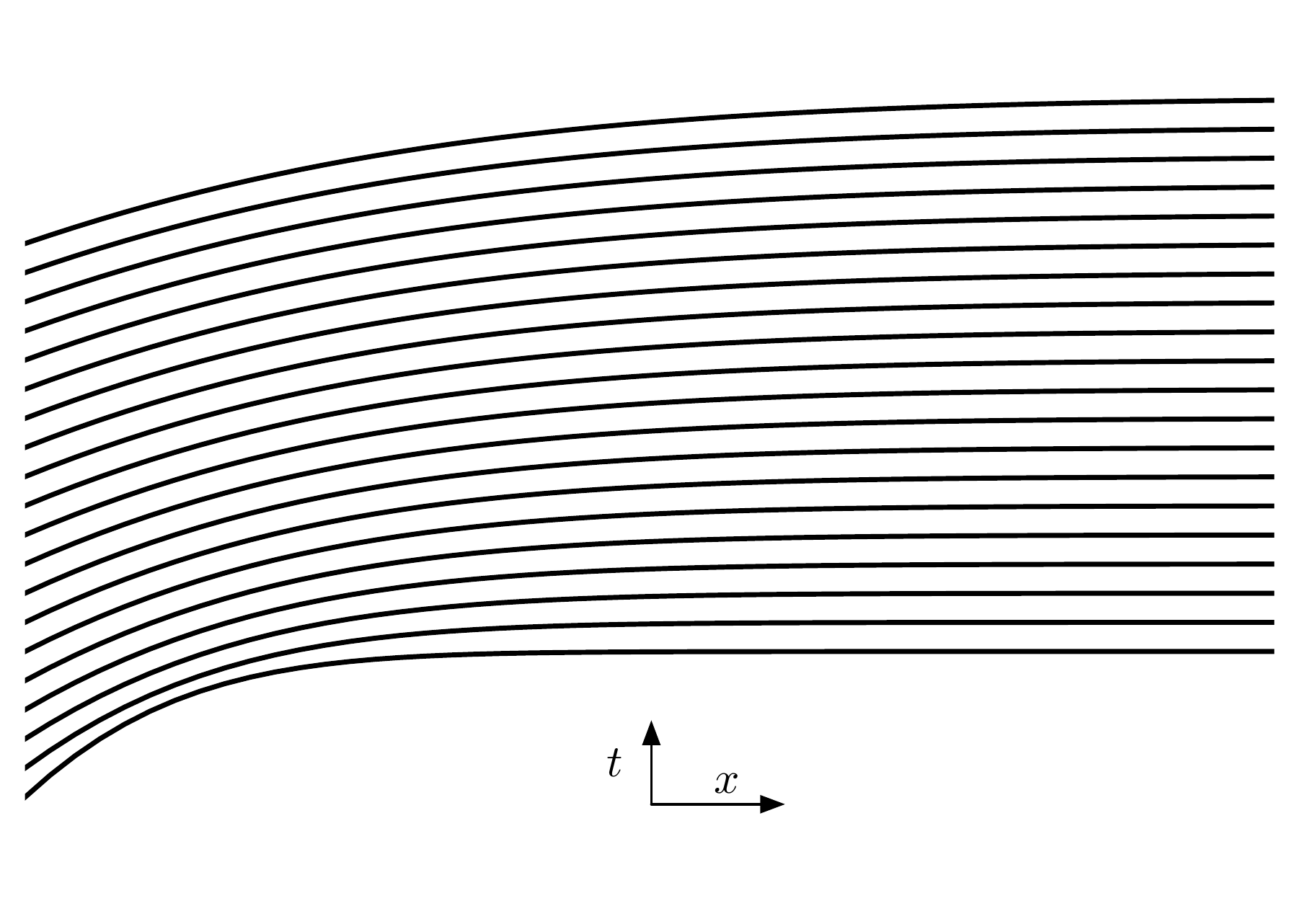}}
\subfigure[]{\includegraphics[width=.49\linewidth]{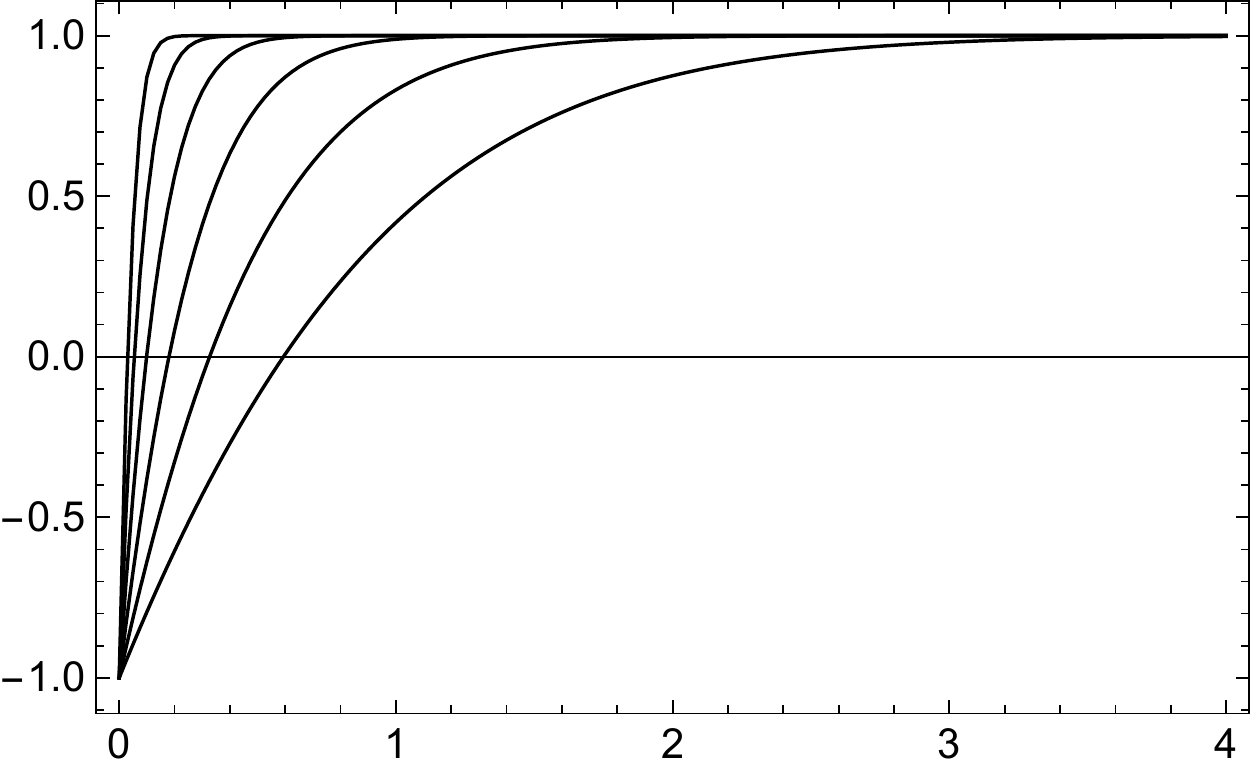}}
\caption{\label{f:Airy1-loc} Plots of $\qloc(x,t)$ for the Airy 1 equation in the concrete case $q_o(0) = 1$ and $g_0(0) = - 1$. (a) The time evolution of $\qloc(x,t)$ up to $t = 2$ for $0 \leq x \leq 15$. (b) An examination of $\qloc(x,t)$ as $t \downarrow 0$ for $t = 1/20(1/6)^j$, $j = 0,1,2,5$.  A discontinuity is formed as $t \downarrow 0$.}
\end{figure}

We construct a similar local solution for \eqref{e:Airy1} where $\omega(k) = -k^3$.  Near a corner singularity we have
\begin{align*}
q(x,t) = 3 q_o(0) I_{\omega,0,1}(x,t) - 3 g_0(0) I_{\omega,0,1}(x,t) + C + \bigo(|x|^{1/2} + |t|^{1/6}).
\end{align*}
To find $C$, we again use that $\lim_{t \downarrow 0}I_{\omega,0,1}(x,t) = 0$ for $x >0$.  Thus $C = q_o(0)$ as above.  We find
\begin{align*}
q(x,t) &= q_{\text{loc}}(x,t) + \bigo(|x|^{1/2} + |t|^{1/6}),\\
\qloc(x,t) &= -3 g_0(0) I_{\omega,0,1}(x,t) + 3q_o(0) \left( I_{\omega,0,1}(x,t)+{1 \over 3} \right).
\end{align*}
We use the same concrete case with the simple data $q_o(0) = 1$ and $g_0(0) = -1$ and we explore $\qloc(x,t)$ in Figure~\ref{f:Airy1-loc}.  Notice that waves travel with a negative velocity because $\omega'(k) < 0$ for $k \in \mathbb R$.  For this reason the corner singularity is regularized for $t\neq 0$ without oscillations.

%%%%%%%%%%%%%%%%%%%%%%%%%%%%%%%%%%%%%%%%%%%%%%%%%%%%%%%%%%%%%%%%%%%%%%%%%%%%%%%%%%%%%%%%%%%%%
\subsection{Airy 2}

 Now, we consider the local solution for \eqref{e:Airy2} where $\omega(k) = k^3$.  Near a corner singularity we have
\begin{align*}
q(x,t) &= q_o(0) \left( (1- \alpha^{2}) I_{\omega,0,2}(x,t) + (1-\alpha) I_{\omega,0,1}(x,t) \right)\\
& + q_0'(0) \left( (1- \alpha) I_{\omega,1,2}(x,t) + (1-\alpha^2) I_{\omega,1,1}(x,t) \right)\\
 &- g_0(0) \left( (1- \alpha^{2}) I_{\omega,0,2}(x,t) + (1-\alpha) I_{\omega,0,1}(x,t) \right)\\
&- g_1(0) \left( (1- \alpha) I_{\omega,1,2}(x,t) + (1-\alpha^2) I_{\omega,1,1}(x,t) \right)\\
&+C_1 + x C_2 + \bigo(|x|^{3/2} + |t|^{1/2}).
\end{align*}
To find $C_1$ we use again use the fact that $\lim_{t \downarrow 0}I_{\omega,i,j}(x,t) = 0$ for $x >0$ and $i \geq 0$.  Thus $C_1 = q_o(0)$.  To find $C_2$ we consider, using \eqref{D1-D2},
\begin{align*}
g_1(0) = \lim_{t \downarrow 0} q_x(0,t) = -3 (g_0(0)-q_o'(0)) I_{\omega,0,1}(0,t) + C_2 + \bigo(|t|^{1/6}).
\end{align*}
But it follows that $I_{\omega,0,1}(0,t) = -1/3$ for $t > 0$ so that $C_2 = q_o'(0)$ and
\begin{align*}
q(x,t) &= \qloc(x,t) + \bigo(|x|^{3/2} + |t|^{1/2}),\\
\qloc(x,t) &= q_o(0)\left( 1+ (1- \alpha^{2}) I_{\omega,0,2}(x,t) + (1-\alpha) I_{\omega,0,1}(x,t) \right)\\
& + q_0'(0) \left(x+ (1- \alpha) I_{\omega,1,2}(x,t) + (1-\alpha^2) I_{\omega,1,1}(x,t) \right)\\
 &- g_0(0) \left( (1- \alpha^{2}) I_{\omega,0,2}(x,t) + (1-\alpha) I_{\omega,0,1}(x,t) \right)\\
&- g_1(0) \left( (1- \alpha) I_{\omega,1,2}(x,t) + (1-\alpha^2) I_{\omega,1,1}(x,t) \right).
\end{align*}

\begin{figure}[htp]
\subfigure[]{\includegraphics[width=.49\linewidth]{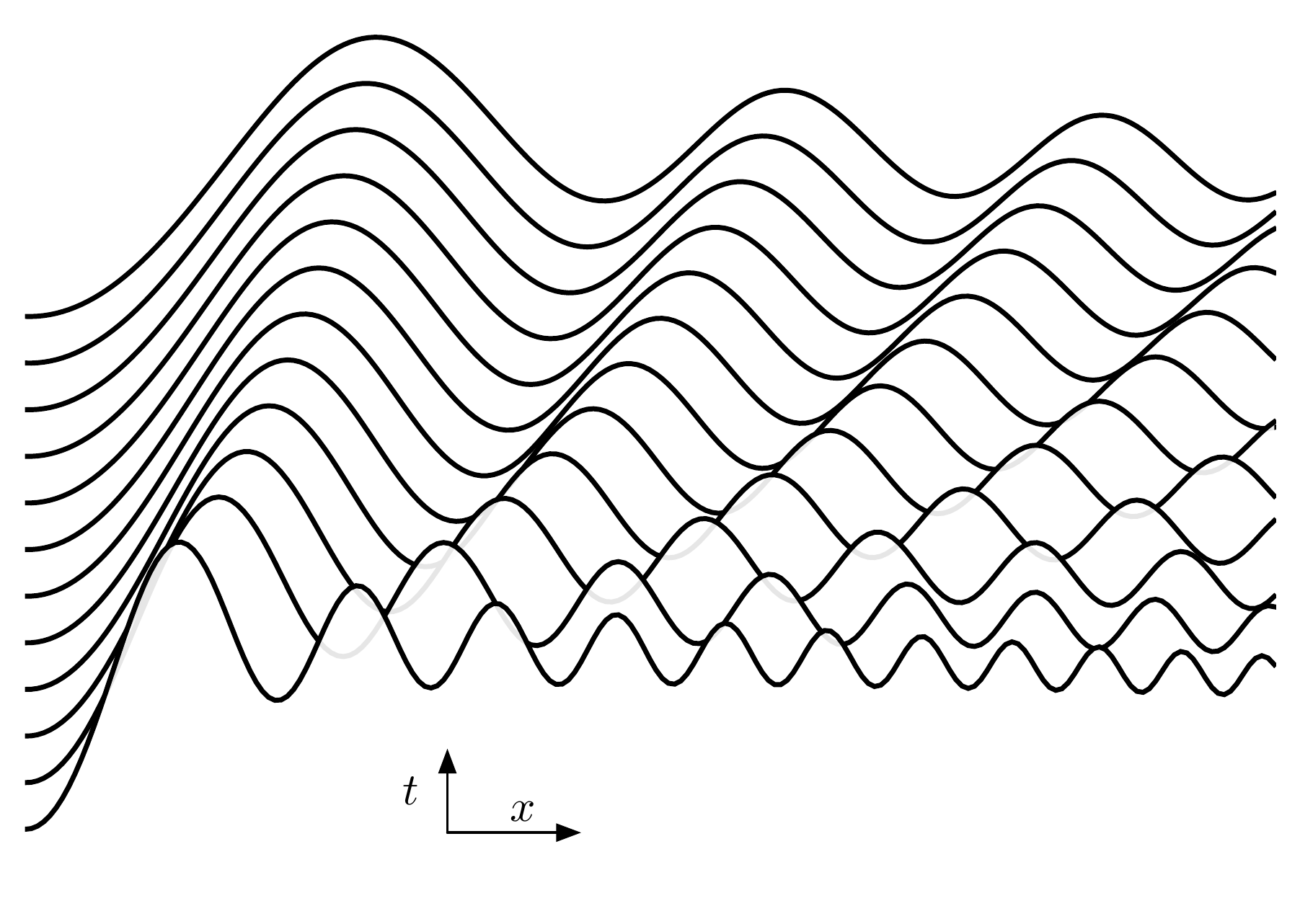}}
\subfigure[]{\includegraphics[width=.49\linewidth]{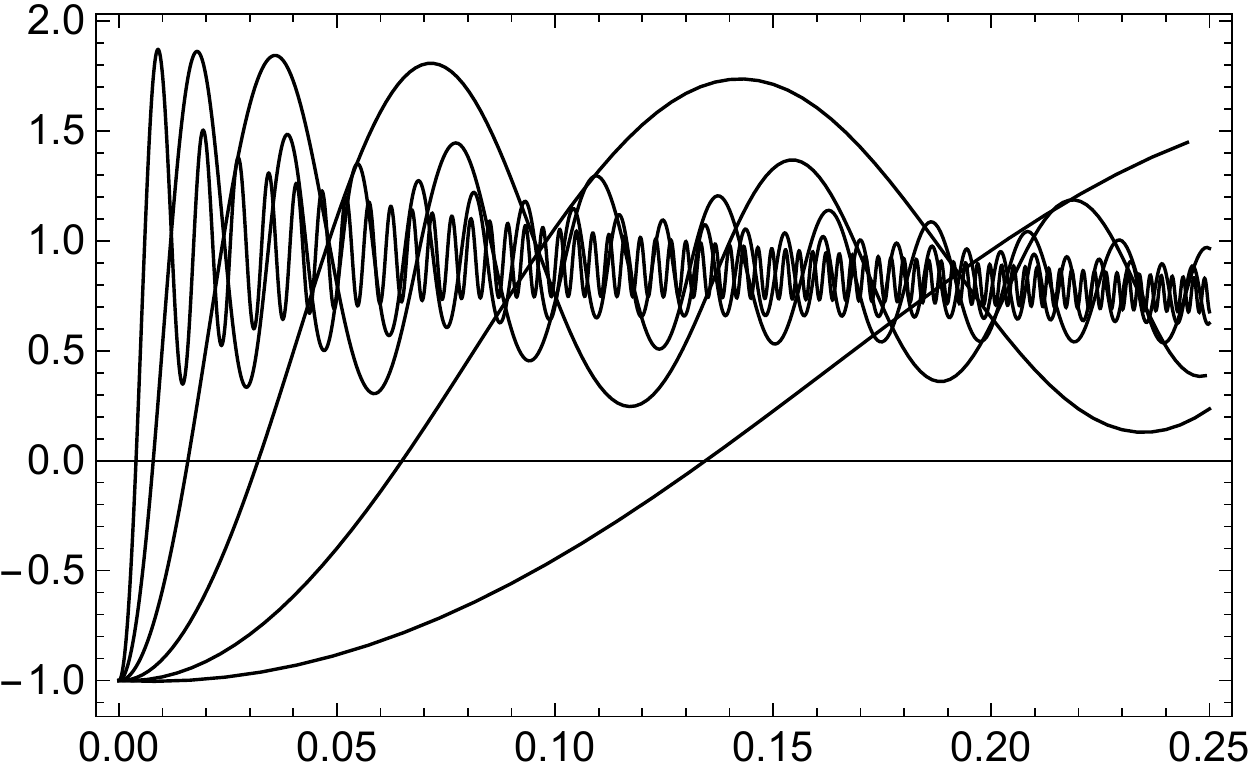}}
\caption{\label{f:Airy2-loc} Plots of $\qloc(x,t)$ for the Airy 2 equation in the concrete case $q_o(0) = 1$, $q_o'(0)=-1$, $g_0(0) = - 1$ and $g_1(0) = 0$. (a) The time evolution of $\qloc(x,t)$ up to $t = 0.00005$ for $0 \leq x \leq 1/2$.  We zoom in on $(x,t) = (0,0)$ in this case so that the effects of the linear term $C_2 x$ are insignificant. (b) An examination of $\qloc(x,t)$ as $t \downarrow 0$ for $t = 1/300(1/8)^j$, $j = 0,1,2,3,4,5$.  A discontinuity is formed as $t \downarrow 0$.}
\end{figure}

\paragraph{First-Order Corner Singularity.} We plot $\qloc(x,t)$ in Figure~\ref{f:Airy2-loc} in the concrete case $q_o(0) = 1$, $q'_o(0) = -1$, $g_0(0) = -1$ and $g_1(0) = -1$.   Note that $q_o'(0) = g_0'(0)$ so that there is no mismatch in the derivative at the origin.

\begin{figure}[htp]
\subfigure[]{\includegraphics[width=.49\linewidth]{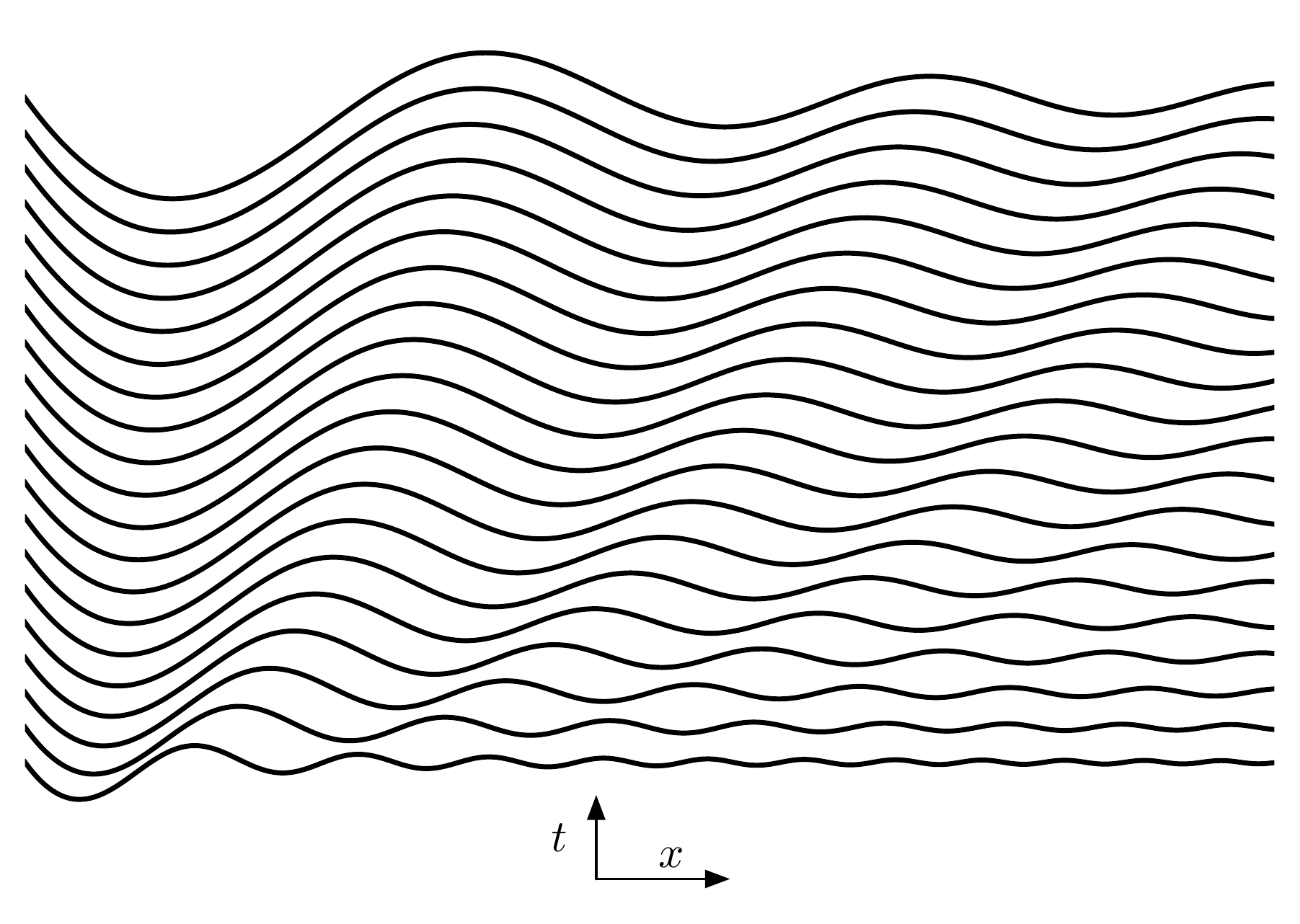}}
\subfigure[]{\includegraphics[width=.49\linewidth]{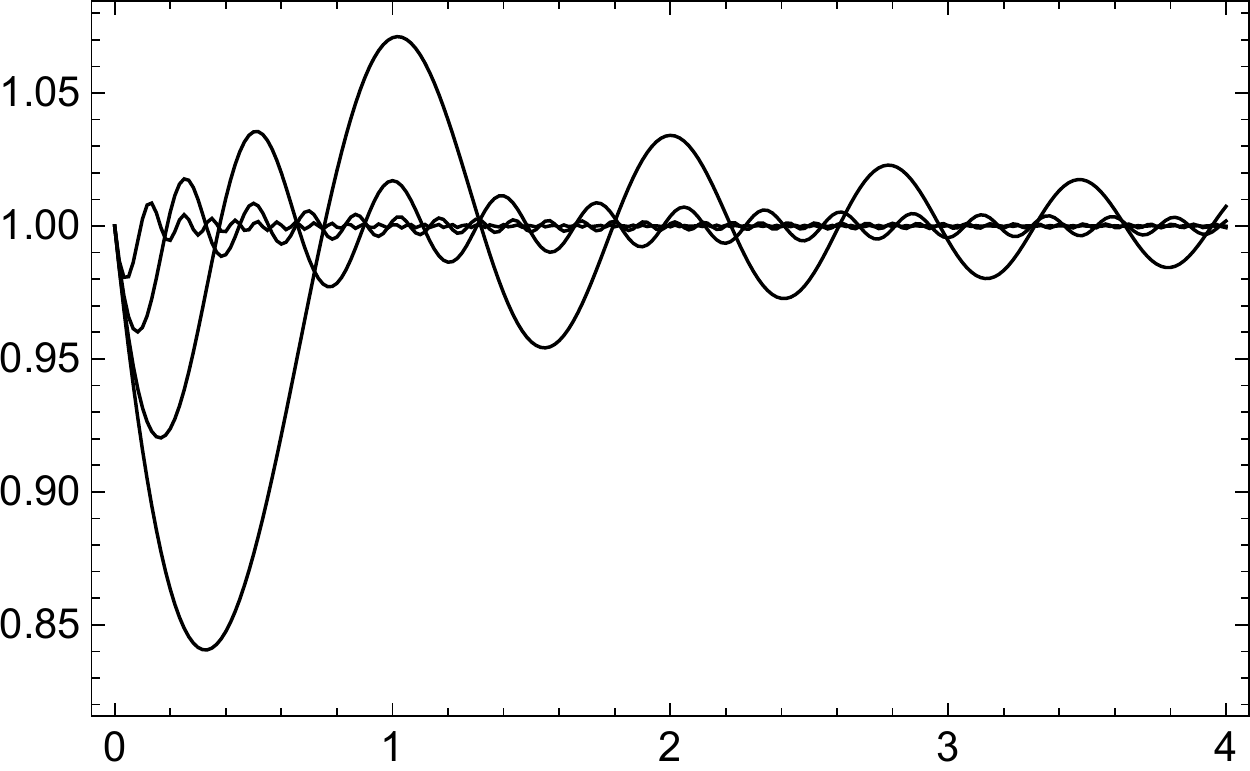}}
\caption{\label{f:Airy2d-loc} Plots of $\qloc(x,t)$ for the Airy 2 equation in the concrete case $q_o(0) = 1$, $q_o'(0)= 0$, $g_0(0) = 1$ and $g_1(0) = -1$. (a) The time evolution of $\qloc(x,t)$ up to $t = 2$ for $0 \leq x \leq 15$. (b) An examination of $\qloc(x,t)$ as $t \downarrow 0$ for $t = 1/10(1/8)^j$, $j = 0,1,2,3,4$.  The function tends uniformly to $q_o(0) = 1$ while $\partial_x q(0,t) = -1$.}
\end{figure}

\paragraph{Second-Order Corner Singularity.} We plot $\qloc(x,t)$ in Figure~\ref{f:Airy2d-loc} in the concrete case $q_o(0) = 1$, $q'_o(0) = 0$, $g_0(0) = 1$ and $g_1(0) = -1$.   Note that $q_o(0) = g_0(0)$ so that there is no mismatch at first order.

\paragraph{An IBVP with discontinuous data.}

\newcommand{\when}{&\text{if}~}
\newcommand{\otherwise}{&\text{otherwise}~}

We now consider the solution of the IBVP for \eqref{e:Airy2} with
\begin{align}\label{disc-data}
\begin{split}
q_o(x) &= \begin{cases} 1, \when x_1 < x < x_2,\\
0, \otherwise,\end{cases}\\
g_0(t) &= \begin{cases} C_1, \when t < t_1,\\
0, \when t \geq t_1,\end{cases}\\
g_1(t) &= C_2.
\end{split}
\end{align}
The solution of this problem has three important features.  The first is the corner singularity at $(x,t) = (0,0)$.  The second is the discontinuities that are present in the initial data.  The last is the singularity in the boundary condition.

\begin{figure}[htp]
\subfigure[]{\includegraphics[width=.49\linewidth]{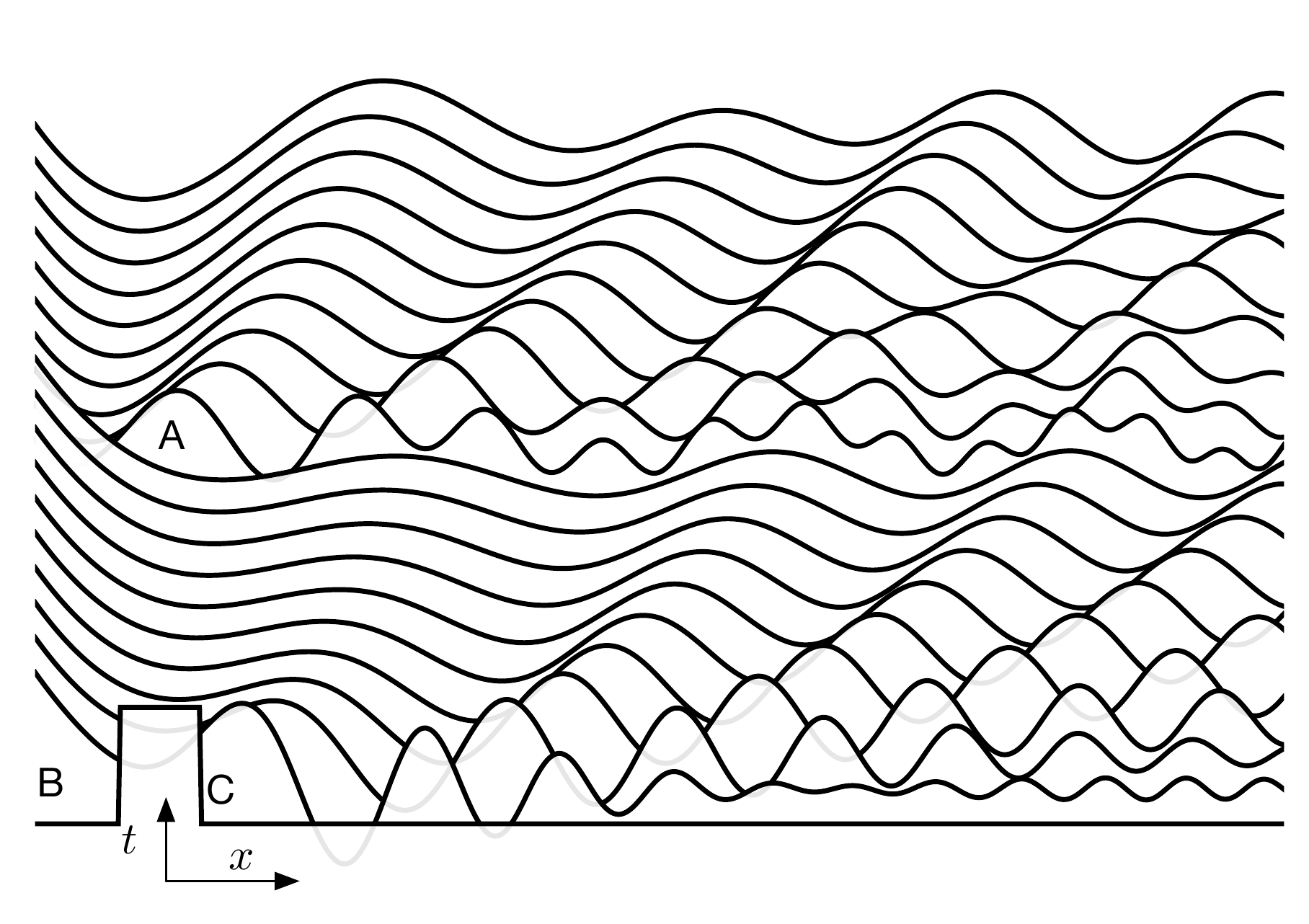}}
\subfigure[]{\includegraphics[width=.49\linewidth]{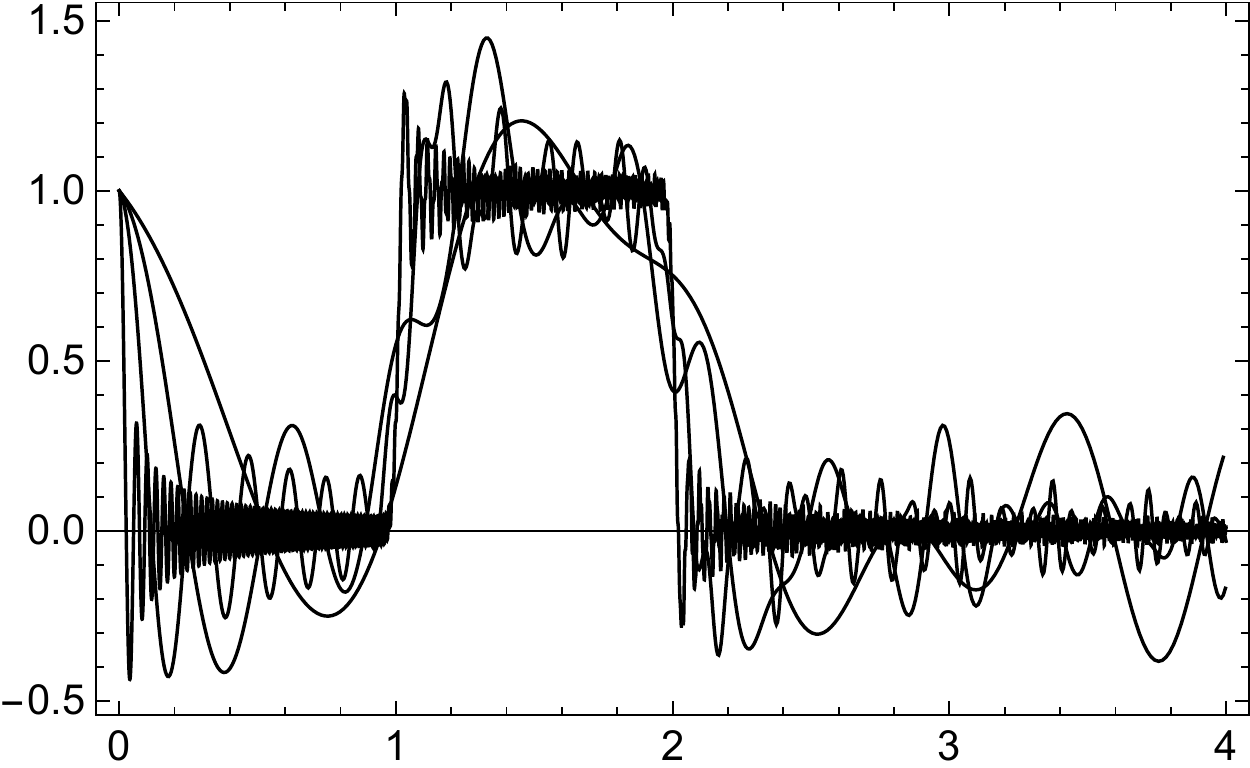}}
\caption{\label{f:Airy2-sol} Plots of $q(x,t)$ for the Airy 2 equation with the data given in \eqref{disc-data} for $C_1 = 1$ and $C_2 = -1$. (a) The time evolution of $q(x,t)$ up to $t = 1$ for $0 \leq x \leq 15$.  Region A signifies the discontinuity in the boundary data, Region B denotes the corner singularity and Region C gives the discontinuity in the initial data. (b) An examination of $q(x,t)$ as $t \downarrow 0$ for $t = 1/10(1/19)^j$, $j = 1,2,3,5$.  A discontinuity is formed as $t \downarrow 0$ at $x = 0,1,2$.}
\end{figure}

Given our developments, this problem can be solved explicitly and computed effectively.  Because $I_{\omega,0,j}(x,t) = 0$ for $t < 0$, the solution formula is
\begin{align*}
q(x,t) &= I_{\omega,0,1}(x-x_1,t) + I_{\omega,0,2}(x-x_1,t) - \alpha^2 I_{\omega,0,2}(x- x_1\alpha,t) - \alpha I_{\omega,0,1}(x-x_1\alpha^2,t)\\
&-I_{\omega,0,1}(x-x_2,t) - I_{\omega,0,2}(x-x_2,t) + \alpha^2 I_{\omega,0,2}(x- x_2 \alpha,t) + \alpha I_{\omega,0,1}(x-x_2 \alpha^2,t)\\
&+C_1 \left( (\alpha^2-1) I_{\omega,0,2}(x,t) + (\alpha-1) I_{\omega,0,1}(x,t) \right)\\
&-C_1 \left( (\alpha^2-1) I_{\omega,0,2}(x,t-t_1) + (\alpha-1) I_{\omega,0,1}(x,t-t_1) \right)\\
&+C_2 \left( (\alpha^2-1) I_{\omega,1,1}(x,t) + (\alpha-1) I_{\omega,1,2}(x,t) \right).
\end{align*}
The solution is plotted in Figure~\ref{f:Airy2-sol}.

\begin{remark}
  For $x > 0$,  $I_{\omega,0,2}(x,t-t_1) = \bigo(|t-t_1|^{1/4})$ as $t \downarrow t_1$ and $I_{\omega,0,2}(x,t-t_1) = 0$ for $t < t_1$.  This implies that $q(x,t)$ is continuous in $t$ but not differentiable at $t = t_1$.  This is a general feature: Discontinuities on the boundary cause the solution to loose time differentiability at that time while the solution maintains continuity.  The above expansions can easily be used to rigorously justify this fact.
\end{remark}

%%% Local Variables:
%%% mode: latex
%%% TeX-master: "gibbs-ibvp.tex"
%%% End:

%%%%%%%%%%%%%%%%%%%%%%%%%%%%%%%%%%%%%%%%%%%%%%%%%%%%%%%%%%%%%%%%%%%%%%%%%%%%%%%%%%%%%%%%%%%%%

%%%%%%%%%%%%%%%%%%%%%%%%%%%%%%%%%%%%%%%%%%%%%%%%%%%%%%%%%%%%%%%%%%%%%%%%%%%%%%%%%%%%%%%%%%%%%
\section*{Acknowledgments}

We thank Bernard Deconinck, Athanassios Fokas, Katie Oliveras, Beatrice Pelloni, Natalie Sheils and Vishal Vasan for many interesting discussions on topics related to this work.  
%We also thank the American Institute of Mathematics in Palo Alto, CA, for its hospitality during the thematic SQuaRE meetings.
This work was partially supported by the National Science Foundation under grant number DMS-1311847 (GB) and DMS-1303018 (TT).

%%%%%%%%%%%%%%%%%%%%%%%%%%%%%%%%%%%%%%%%%%%%%%%%%%%%%%%%%%%%%%%%%%%%%%%%%%%%%%%%%%%%%%%%%%%%%
\begingroup
\renewcommand\thesubsection{A.\Roman{subsection}}
\renewcommand\theequation{\Alph{subsection}.\arabic{equation}}
\setcounter{section}0
\setcounter{subsection}0
\setcounter{equation}0
\def\thetheorem{A.\arabic{theorem}}

\section*{Appendix}

\subsection{Validity of solution formula and regularity results}
\label{a:SolutionForm}

From the work of \cite{FokasSung} we know that the expression \eqref{e:form} evaluates to give the solution of \eqref{generalPDE} pointwise provided the initial and boundary data are sufficiently regular.

\begin{lemma}\label{Lemma:IntOrder}
If $g_j \in H^1(0,T)$ and $q_0 \in L^1 \cap L^2(\mathbb R)$ each integral in \eref{e:form} can be written in the form
\begin{align*}
g_j(t) T_j(x,t,t) - g_j(0) T_j(x,t,0) - \int_{0}^t T_j(x,t,s) g'_j(s) ds, ~~\text{or}~~ \int_{0}^\infty S(x,t,s) q_0(s) ds,
\end{align*}
where $S(x,t,s)$ and $T_j(x,t,s)$ are bounded in $s$ for fixed $x>0$ and $t>0$.  \br Furthermore, for $\kappa = 0,1,2,\ldots$
\begin{itemize}
\item $\partial_x^\kappa S(x,t,s) \sim |s|^{\frac{2\kappa-n+2}{2(n-1)}}$ as $s \rightarrow \infty$,
\item $\partial_t^\kappa S(x,t,s) \sim |s|^{\frac{2n\kappa-n+2}{2(n-1)}}$ as $s \rightarrow \infty$,
\item $\partial_x^\kappa T(x,t,s) \sim |s-t|^{\frac{n+2j-2\kappa}{2(n-1)}}$ as $s \rightarrow t^-$, and
\item $\partial_t^\kappa T(x,t,s) \sim |s-t|^{\frac{n+2j-2n\kappa}{2(n-1)}}$ as $s \rightarrow t^-$.
\end{itemize}
\er

\begin{proof}
The estimate for the integral
\begin{align*}
\frac{1}{2\pi} \int_{\mathbb R} e^{ik(x-s)-i\omega(k)t} dk
\end{align*}
which is the kernel in the integral
\begin{align*}
\frac{1}{2 \pi} \int_{\mathbb R} e^{ikx - i\omega(k)t} \hat q_o(k) dk
\end{align*}
follows directly from Theorem~\ref{Thm:Kernel}.  Next consider the integral 
\begin{align*}
\int_{\partial D_i^+} e^{ikx-i\omega(k)t} \hat q(\nu(k)) dk &= \lim_{R \rightarrow \infty} \int_{\partial D_i^+ \cap B(0,R)} e^{ikx-i\omega(k)t} \hat q(\nu(k)) dk\\
&= \lim_{R \rightarrow \infty} \int_{0}^\infty S_R(x,t,s) q_0(s) ds, \\
S_{R}(x,t,s) &= \int_{\partial D_i^+ \cap B(0,R)} e^{ikx-i\nu(k)s-i\omega(k)t}dk.
\end{align*}
We perform a change of variables on $S_R$
\begin{align*}
S_{R}(x,t,s) &= \int_{\nu^{-1}(\partial D_i^+ \cap B(0,R))} e^{-izs + i\nu^{-1}(z) x-i\omega(z)t} d \nu^{-1}(z).
\end{align*}
Here $\nu^{-1}(D_i^+)$ is a component of $D$ in $\mathbb C^-$.  We discuss the case where $\nu(k) = \alpha k$ for $|\alpha| = 1$, \emph{i.e.} $\omega(k) = \omega_n k^n$.  The general case follows from similar but more technical arguments.  For fixed $x$ and $t$ we apply Theorem~\ref{Thm:Kernel} with $w(k) = \omega(z) - \alpha^{-1} z x/t$ after possible deformations.  In all cases, $e^{-izs + i\nu^{-1}(z) x-i\omega(z)t}$ is bounded large $s$ when $z$ is replaced with the appropriate stationary phase point. We obtain
\begin{align*}
\lim_{R \rightarrow \infty} \partial_x^j S_R(x,t,s) \sim |s|^{\frac{2j +2-n}{2(n-1)}}.
\end{align*}

Next, we consider the terms involving $g_j$.  Generally speaking, for the canonical problem with $\omega(k) = \omega_n k^n$ these terms are of the form
\begin{align*}
\int_{\partial D_i^+} e^{ikx-i\omega(k)t} k^{N(n)-j} \tilde g_j(-\omega(k),t) dk  &= \lim_{R \rightarrow \infty} \int_{\partial D_i^+ \cap B(0,R)} e^{ikx-i\omega(k)t} k^{N(n)-j} \tilde g_j(-\omega(k),t) dk.
\end{align*}
We write
\begin{align*}
e^{ikx-i\omega(k)t} k^{N(n)-j} \tilde g_j(-\omega(k),t) = e^{ikx-i\omega(k)t} \frac{k^{N(n)-j}}{i\omega(k)} \left (  g_j(t)e^{i \omega(k) t} -g_j(0)  - \int_0^t e^{i \omega(k) s} g_j'(s)ds \right)
\end{align*}
so that
\begin{align*}
\int_{\partial D_i^+ \cap B(0,R)} e^{ikx-i\omega(k)t} k^{N(n)-j} \tilde g_j(-\omega(k),t) dk &= \frac{g_j(t)}{i\omega_n} \int_{\partial D_i^+ \cap B(0,R)} e^{ikx} \frac{dk}{k^{n-N(n)+j}} \\
- \frac{g_j(0)}{i\omega_n} \int_{\partial D_i^+ \cap B(0,R)} e^{ikx-i\omega(k)t} \frac{dk}{k^{n-N(n)+j}}
&- \int_{0}^t \left( \frac{1}{i\omega_n} \int_{\partial D_i^+ \cap B(0,R)} e^{ikx-i\omega(k)(t-s)} \frac{dk}{k^{n-N(n)+j}} \right) g_j'(s) ds.
\end{align*}
Now, because $n-N(n)+j \geq 1$ all integrals converge for $x > 0$ as $R \rightarrow \infty$. Additionally, the integral with $g_j(t)$ as a coefficient vanishes identically. For $x > 0$ by Theorem~\ref{Thm:Kernel} \br with $m = n - N(n) + j -1$
\begin{align*}
\lim_{R \rightarrow \infty} \int_{\partial D_i^+ \cap B(0,R)} e^{ikx-i\omega(k)(t-s)} \frac{dk}{k^{n-N(n)+j}} = \mathcal O \left(  |s-t|^{\frac{n + 2(j-1)}{2(n-1)}}\right),
\end{align*}\er
as $s \rightarrow t^-$, implying this is a bounded function for all $s \in [0,t]$. To estimate $t$ derivatives we note that the estimates for $\partial^{jn}_x$ follow for $\partial_t^j$. This proves the lemma.

\end{proof}
\end{lemma}

\begin{lemma}\label{Lemma:Valid}
The solution formula~\eqref{e:form} holds for $q_0 \in L^1 \cap L^2(\mathbb R^+)$ and $g_j \in H^{1}([0,T])$ for all $t > 0, x >0$, $j = 0, \ldots, N(n)-1$.
\begin{proof}
To prove this result we must approximate $q_0$ and $g_j$ with smooth functions that are compatible at $(x,t) = (0,0)$. First, we find a sequence of functions $\tilde q_{0,n} \in C_c^\infty((0,R))$ such that $q_{0,n} \rightarrow q_{0}$ in $L^1 \cap L^2(\mathbb R^+)$ .  To see that such a sequence exists, consider the approximation of $q_0(x) \chi_{[0,R]}(x)$ in $L^2(\mathbb R^+)$ with $C_c^\infty((0,R))$ functions.  Because of the bounded interval of support, this approximation converges in $L^1(\mathbb R^+)$ as well.  Next, because $q_0(x) \chi_{[0,R]}(x) \rightarrow q_0(x)$ in $L^1\cap L^2(\mathbb R^+)$ as $R \rightarrow \infty$, a diagonal argument produces an acceptable sequence.   Now, find sequences $d_{j,n} \rightarrow g_j'$ in $L^2(0,T)$ with $d_{j,n} \in C_{c}^\infty(0,T)$. Then define
\begin{align*}
g_{j,n}(t) = g_{j}(0) + \int_0^t d_{j,n}(s) ds,
\end{align*}
so that $g_{j,n}$ is constant near $t = 0$.  Define $p(x) = \sum_{j=0}^{N(n)-1} g_j(0) \frac{x^j}{j!}$ and $\phi_{n}(x)$  have support $[0,2/n]$ and be equal to $1$ on $[0,1/n]$ and interpolate smoothly and monotonically between $0$ and $1$ on $[1/n,2/n]$.  Then $q_{0,n}(x) + p(x)\phi_n(x)$ converges to $q_0$ in $L^2(\mathbb R^+)$ and $q_{0,n}$ and $g_{j,n}$ are compatible at $(x,t) = (0,0)$ and the solution formula \eref{e:form} holds with this combination of initial/boundary data.

Now, because convergence of the initial data also occurs in $L^1(\mathbb R^+)$ and convergence of the boundary data also occurs in\footnote{$W^{1,1}(0,T)$ is the space of integrable functions on the interval $(0,T)$ with one integrable derivative.} $W^{1,1}(0,T)$, we apply Lemma~\ref{Lemma:IntOrder} to demonstrate that the solution formula with data $(q_{0,n},g_{j,n})$ converges pointwise to the solution value and furthermore limits may be passed inside the relevant integrals.  This implies the solution formula holds with these relaxed assumptions.
\end{proof}
\end{lemma}

To handle multiple boundary discontinuities, we note that we can solve the problem with zero initial data. Assume the boundary condition has a discontinuity at $0 < t_1 < T$.  With boundary conditions 
\begin{align*}
g_j(t) = \begin{cases} g_{j,1}(t),& t \in [0,t_1],\\
g_{j,2}(t), & t \in (t_1,T].
\end{cases}
\end{align*}
That are piecewise $H^1$ functions.  We use linearity to modify the boundary condition.  Consider the two functions
\begin{align*}
G_{j,1}(t) &= \begin{cases} g_{j,1}(t), & t \in [0,t_1],\\
g_{j,1}(t_1),& t \in (t_1,T],
\end{cases}\\
G_{j,2}(t) &= \begin{cases} 0, & t \in [0,t_1],\\
g_{j,2}(t)-g_{j,1}(t_1),& t \in (t_1,T],
\end{cases}\\
\end{align*}

Since the above theorem indicates the solution is given by the formula for all $t \in [0,T]$, with boundary conditions $G_{j,1}$.   Furthermore, the initial-boundary-value problem with zero initial data and boundary data $G_{j,2}$ is also given by the solution formula, with the solution being identically zero before $t = t_1$.  We use linearity to add these two solutions. We have shown that the \eref{e:form} gives us this weak solution in the interior.

Further considerations can be used to show the solution is smooth in $x$ for all $t > 0$ and smooth in $t$ for $t> 0$, $t \neq t_1$.  The contributions from integrals involving $g_j$ can cause complicated singularities in the solution.  With this in mind we state our regularity theorem.

\begin{theorem}\label{Thm:Reg}
Assume $q_0 \in L^2(\mathbb R^+) \cap L^1(\mathbb R^+,(1+|x|)^\ell)$ and $g_j \in H^{p+1}(t_i,t_{i+1})$ ($p \geq 0$) for $0 = t_0 < \cdots < t_m = T$.  Then \eref{e:form} evaluates pointwise to give the $L^2$ solution of \eref{generalPDE}.
\begin{itemize}
\item If
\begin{align*}
\ell \geq \frac{2m -n + 2}{2(n-1)}, \quad np - [2 N(n) + 2 -n] >  m,
\end{align*}
then $q(x,t)$ is differentiable $m$ times with respect to $x$ for $x> 0, t >0$.
\item If 
\begin{align*}
\ell \geq \frac{2jn -n +2}{2(n-1)}, \quad p - \frac{1}{n}[2 N(n) + 2 -n] >  j
\end{align*}
then $q(x,t))$ is differentiable $j$ times with respect to $t$ for $x > 0, ~t \neq t_i$ and continuous in $t$ for $t > 0$.
\end{itemize}
\begin{proof}
Lemma~\ref{Lemma:Valid} demonstrates that \eref{e:form} produces the solution pointwise for $t \leq t_1$.  We look at the differentiability of the solution in $(0,\infty) \times (0,t_1)$.  The differentiability of the integrals in \eqref{e:form} that involve $q_0(x)$ follows from the growth of the kernel.  To see differentiability in of the terms involving $g_j$ we note that integration by parts can be performed $p$ times.  It remains to consider the differentiability.  \br Formally,
\begin{align*}
  \frac{d^\kappa}{dx^\kappa} \int_{0}^t \left(\int_{\partial D_i^+} e^{ikx - i\omega(k)(t-s)} \frac{dk}{k^{pn-N(n)+j}} \right) g_j^{(p+1)}(s) ds \\
  = i^\kappa \int_{0}^t \left(\int_{\partial D_i^+} e^{ikx - i\omega(k)(t-s)} \frac{dk}{k^{pn-N(n)+j-\kappa}} \right) g_j^{(p+1)}(s) ds.
\end{align*}
It is straightforward to check from Theorem~\ref{Thm:Kernel} that the kernel in this integral is an $L^2$ function, and differentiability follows, provided $2pn - 2 N(n) + 2j -2\kappa -1 >  1 - n$ and for simplicity a condition is $\kappa < pn - [2 N(n) + 2 -n]$. \er   This implies we may take $pn$ $x$-derivatives inside the integral and $p$ $t$-derivatives.

Next define $G_{j,1}(t)$ to be an $H^{p+1}((0,T))$ extension of $g_j(t)\chi_{[0,t_1]}(t)$.  Iteratively, define
\begin{align*}
G_{j,i}(t) = g_j(t) - \sum_{M=1}^{i-1} G_{j,M}(t), \quad t \in [t_{i-1},t_i), \quad i = 2, \ldots, m,
\end{align*}
and assume each of these are extended as an $H^{p+1}((t_{i-1},T))$ function.  Let $q_i(x,t)$ be the solution of \eref{generalPDE} with initial/boundary data $(q_0,G_{j,1})$ if  $i = 1$ and $(q_{0,i} \equiv 0, G_{j,i})$ for $i > 1$ on $(0,\infty) \times(t_i,T]$. The solution formula \eref{e:form} is valid with this initial data.  The solution with data $(q_0,g_j)$ is given by
\begin{align*}
q(x,t) = \sum_{M=1}^i q_i(x,t), \quad t \in [t_{i-1},t_i),
\end{align*}
and the regularity follows.
\end{proof}
\end{theorem}
\begin{remark}
  This theorem can be improved by allowing $p = p(j)$ to be a fraction.  But our aim is only to give sufficient conditions for differentiability that are simple to state.
\end{remark}

%%% Local Variables:
%%% mode: latex
%%% TeX-master: "gibbs-ibvp"
%%% End:

\subsection{Special functions arising in the IBVP}\label{a:SF}

Recall
\begin{align*}
I_{w,m,j}(x,t) = \frac{1}{2\pi} \int_{\partial D^+_j} \frac{ e^{ikx-iw(k)t}}{(ik)^{m+1}} dk,
\end{align*}
and suppose $w(k) = w_n k^n + \bigo( k^{n-1})$.  Further, define
\begin{align*}
K_t(x) = \sum_{j=1}^{N(n)} I_{w,-1}(x,t).
\end{align*}
For $|x| > 0$, $t > 0$, we rescale, by setting $\sigma = \sign(x)$, $k = \sigma (|x|/t)^{1/(n-1)}z$
\begin{align}\label{e:rescale}
\begin{split}
I_{\omega,m,j}(x,t) &= \sigma^m \left(\frac{|x|}{t}\right)^{-m/(n-1)} \int_{\Gamma_j} e^{X(iz-i\omega_n \sigma^n z^n-i R_{|x|/t}(z))} \frac{dz}{(iz)^{m+1}},\\
R_{|x|/t}(z) &= \sum_{j=2}^{n-1} \omega_j \left(\frac{|x|}{t}\right)^{\frac{j-n}{n-1}} (\sigma z)^j, \quad X = |x| \left( \frac{|x|}{t} \right)^{1/(n-1)}.
\end{split}
\end{align}
Define
\begin{align*}
\Phi_{|x|/t}(z) = ik-i\omega_n \sigma^n z^n-i R_{|x|/t}(z),
\end{align*}
where $\{z_j\}_{j=1}^{n-1}$ are the roots of $\Phi_{|x|/t}'(z) = 0$ ordered counter-clockwise from the real axis. Here $\Gamma_j$ is a deformation of $\partial D_i^+$ which passes along the path of steepest descent through $z_j$.

\begin{theorem}\label{Thm:Kernel}
Suppose $\omega(k) = \omega_n k^n + \bigo( k^{n-1})$ then as $|x/t| \rightarrow \infty$
\begin{align*}
I_{\omega,m,j}(x,t) &= -i \Res_{k=0}  \left(\frac{ e^{ikx-i\omega(k)t}}{(ik)^{m+1}}\right) \chi_{(-\infty,0)}(x) \\
&+  \frac{\sigma^m |x|^{-1/2}}{\sqrt{2 \pi}}\left(\frac{|x|}{t}\right)^{-\frac{m+1/2}{n-1}}   \frac{e^{X \Phi_{|x|/t}(z_j) +i\theta_j}}{(i z_j)^{m+1}} \frac{1}{|\Phi_{|x|/t}''(z_j)|^{1/2}}\left(1 + \bigo\left(|x|^{-1} \left( \frac{|x|}{t} \right)^{-1/(n-1)}\right)\right).
\end{align*}
Here $\theta_j$ is the direction at which $\Gamma_j$ leaves $z_j$.  Hence
\begin{itemize}
\item  \br For fixed $t> 0$ as  $|x| \rightarrow \infty$
\begin{align}\label{e:largespace}
K^{(m)}_t(x) = \left\{ \begin{array}{lr} \bigo \left(|x|^{\frac{2m-n+2}{2(n-1)}} \right), & n \text{ is even},\\
\\
\bigo \left(|x|^{\frac{2m-n+2}{2(n-1)}}\right), & n \text{ is odd}, ~~ \omega_n x > 0,\\
\\
\bigo \left(|x|^{-M} \right) \text{ for all } M >0, & n \text{ is odd}, ~~ \omega_n x < 0.
\end{array}\right.
\end{align}\er
\item  For $|x| \geq \delta > 0$ and $m \geq 0$ as $t \rightarrow 0^+$
\begin{align}\label{e:smalltime}
I_{\omega,m}(x,t) &= -i \Res_{k=0}  \left(\frac{ e^{ikx-i\omega(k)t}}{(ik)^{m+1}}\right) \chi_{(-\infty,0)}(x) + \bigo\left(t^{\frac{m+1/2}{n-1}} |x|^{-\frac{2m+2n}{2(n-1)}} \right).
\end{align}

\end{itemize}
\end{theorem}

%%% Local Variables:
%%% mode: latex
%%% TeX-master: "gibbs-ibvp"
%%% End:

\subsection{Residual estimation}
\label{a:Residual}

In many cases we must understand the behavior of integrals of the form
\begin{align*}
\int_S e^{ikx-i\omega(k)t} F(k) \frac{dk}{k^{m+1}}
\end{align*}
for small $|x|$ and $t$.  Here $S$ is a piecewise smooth, asymptotically affine contour in the upper-half plane that avoids the origin along which $e^{-i\omega(k)t}$ is bounded.  One might expect that a Taylor expansion of the integrand near zero would provide the leading contribution. Namely,
\begin{align*}
\int_S e^{iks-i\omega(k)t} F(k) \frac{dk}{k^{m+1}} = \sum_{j=0}^m \int_{S} a_j(x,t) k^{j-m-1} F(k)dk + E_m(x,t),
\end{align*}
where $E_m$ is of higher-order as $(x,t) \rightarrow (0,0)$.  We make this fact rigorous in this section.  Define $a_j(x,t)$ to be the $j$th-order Taylor coefficient of $e^{ikx-i\omega(k)t}$ at $ k = 0$. We make some observations about these coefficients.  We write
\begin{align*}
ikx - i\omega(k)t = ikx -i \sum_{j=2}^n \omega_j (t^{1/j} k)^j.
\end{align*}
From this it is clear that $|a_j(x,t)| \leq C \sum_{p=0}^j |x|^{p} t^{\frac{j-p}{n}}$.  With each power of $k$ comes a power of $x$ or a least $t^{1/n}$.  Define $\rho(x,t) = |x| + |t|^{1/n}$ and there exists $C_j > 0$ such thaty
\begin{align}\label{e:TaylorTerm}
\frac{1}{C_j} \rho(x,t)^j \leq \sum_{p=0}^j |x|^{p} t^{\frac{j-p}{n}} \leq C_j \rho(x,t)^{j}.
\end{align}
We also want to understand the behavior of the derivatives of $e^{ikx-i\omega(k)t}$ in the complex plane.  Namely, we want to understand which powers of $x$ and $t$ go with powers of $k$.  The first few derivatives are, of course,
\begin{align*}
(ix &- i \omega'(k) t) e^{ikx-i\omega(k)t},\\
(ix &- i \omega'(k) t)^2e^{ikx-i\omega(k)t}  + (-i \omega''(k) t) e^{ikx-i\omega(k)t}.\\
(ix &- i \omega'(k) t)^3 e^{ikx-i\omega(k)t} + 2(-i \omega''(k) t) e^{ikx-i\omega(k)t} + (-2 i \omega''(k)t) e^{ikx-i\omega(k)t}.
\end{align*}
The observation to be made here is that for $|k| \geq 1$, $|x|,t \leq 1$ there are positive constants $D_j$ and $B_j$ such that
\begin{align}
\left| \frac{d^j}{dk^j} e^{ikx-i\omega(k)t} \right| &\leq D_j \left( |x| + n t \sum_{p=2}^{n} |\omega_n| |k|^{p-1} \right)^j \left|e^{ikx-i\omega(k)t}\right|\notag\\
& \leq B_j \rho(x,t)^j (1+ \rho(x,t)|k|)^{j(n-1)}\left|e^{ikx-i\omega(k)t}\right|.\label{e:TaylorError}
\end{align}
These are the necessary components to prove the following.

\begin{lemma}\label{Lemma:Expansion}
Suppose $S$ be a piecewise smooth, asymptotically affine contour in the upper-half plane, avoiding the origin, such that $e^{-i\omega(k)t}$ is bounded on $S$ for $0 \leq t \leq 1$.  If $F \in L^2(S)$ there exists a constant $C > 0$ such that
\begin{align*}
\left| \int_S e^{ikx-i\omega(k)t} F(k) \frac{dk}{k^{m+1}} - \sum_{j=0}^m \int_{S} a_j(x,t) k^{j-m-1} F(k)dk \right| \leq C\rho^{m+1/2}(x,t) \|F\|_{L^2(S)}.
\end{align*}

\begin{proof}
Define
\begin{align*}
f_{x,t,m}(k) = \frac{1}{k^{m+1}} \left(  e^{ikx-i\omega(k)t} - \sum_{j=0}^m  a_j(x,t) k^{j}  \right).
\end{align*}
We estimate the $L^2(S)$ norm of this function.  First for $\rho \equiv \rho(x,t)$, $k \in S \cap B(0, \rho^{-1})$ we have by Taylor's Theorem applied along $S$ (using its smoothness) there exists $C_m > 0$ such that (see \eref{e:TaylorError})
\begin{align*}
\left| e^{ikx-i\omega(k)t} - \sum_{j=0}^m a_j(x,t) k^{j}\right| \leq C_m\frac{|k|^{m+1}}{(m+1)!} \rho^{m+1} \sup_{k\in S} \left|e^{ikx-i\omega(k)t}\right|.
\end{align*}
From this we find that for a (new) constant $C_m > 0$
\begin{align}\label{e:inner-est}
\left( \int_{S \cap B(0,\rho^{-1})} |f_{x,t,m}(k)|^2 |dk| \right)^{1/2} \leq \frac{C_m}{(m+1)!} \rho^{m+1/2},
\end{align}
because $ \int_{S \cap B(0,R)}|dk| = \bigo(R)$ as $R \rightarrow \infty$.

Next, we estimate on $S \setminus B(0,\rho^{-1})$.  In general, we find
\begin{align*}
\left( \int_{S \setminus B(0,\rho^{-1})} |k|^{2(j-m-1)} |dk| \right)^{1/2} \leq D_j \rho^{m-j+1/2},
\end{align*}
and using \eref{e:TaylorTerm}
\begin{align}\label{e:outer-est}
\left( \int_{S \setminus B(0,R)} |f_{x,t,m}(k)|^2 |dk|\right)^{1/2} \leq C \sum_{j=0}^\infty D_j \rho^{m+1/2}.
\end{align}
Combining \eref{e:inner-est} and \eref{e:outer-est} with the Cauchy--Schwarz inequality proves the result.

\end{proof}
\end{lemma}

The final piece we need is sufficient conditions for $F \in L^2(S)$.  Recall that $S$ is always in the domain of analyticity of
\begin{align*}
F(\nu(k)) = \int_0^\infty e^{-i\nu(k)x} f(x) dx.
\end{align*}
More precisely, $\nu^{-1}(S)$ is in the closed lower-half plane.  So
\begin{align*}
\int_{S} |F(\nu(k))|^2 |dk| = \int_{\nu^{-1}(S)} |F(k)|^2 |d\nu^{-1}(k)|.
\end{align*}
Also, $S$ can be chosen such that $\nu^{-1}$ has a uniformly bounded derivative on $\nu^{-1}(S)$ (see \cite{DTV}).  It follows that $F$ is in the Hardy space of the lower-half place (see \cite{TrogdonThesis}) and can be represented as the Cauchy integral of its boundary values
\begin{align*}
\mathcal C_{\mathbb R} F(k) = \frac{1}{2 \pi i} \int_{\mathbb R} \frac{F(z)}{z-k} dz = - F(k).
\end{align*}
The Cauchy integral operator is bounded on $L^2(\mathbb R \cup S)$ so that
\begin{align*}
\|F\|_{L^2(S)}= \|\mathcal C_{\mathbb R}F\|_{L^2(S)} \leq \|\mathcal C_{\mathbb R}F\|_{L^2(\mathbb R \cup S)} \leq C \|F\|_{L^2(\mathbb R)}.
\end{align*}

Next, $S$ is always in the domain of analyticity and boundedness of
\begin{align*}
G(-\omega(k) ) =\int_{0}^t e^{i\omega(k)s}g(s) ds.
\end{align*}
This is true because $S$ asymptotically is a subset of $\partial D_i^+$.  Set $z = -\omega(k)$, noting that $z \in \mathbb C^-$, we have
\begin{align*}
\int_{S} |G(-\omega(k)|^2 |d(\omega(k))| = \int_{-\omega(S)} |G(z)|^2 dz < \infty,
\end{align*}
if $g \in L^2(0,t)$.   Furthermore, if $S$ avoids zeros of $\omega'$
\begin{align*}
\int_{S} |G(-\omega(k)|^2 |dk| \leq C' \int_{S} |G(-\omega(k)|^2 |d(\omega(k))|, \quad C' > 0.
\end{align*}
Similar Hardy space considerations indicate that if $g \in L^2(0,t)$ then $G(-\omega(\cdot)) \in L^2(S)$.  We obtain the following\footnote{Such a theorem holds on contours with much less regularity.}.
\begin{lemma}\label{Lemma:Sufficient}
\begin{itemize}
Let $S$ be a Lipschitz contour.
\item If $f \in L^2(\mathbb R^+)$, $\im \nu(k) \leq 0$ on $S$ and  $\nu^{-1}$ has a uniformly bounded derivative on $\nu(S)$, then $F \in L^2(S)$.
\item If $g \in L^2(0,t)$ and $S \subset D$ is bounded away from the zeros of $\omega'$, then $G(-\omega(k)) \in L^2(S, |d(\omega(k))|) \subset L^2(S)$.

\end{itemize}
\end{lemma}

%%% Local Variables:
%%% mode: latex
%%% TeX-master: "gibbs-ibvp.tex"
%%% End:

%%%%%%%%%%%%%%%%%%%%%%%%%%%%%%%%%%%%%%%%%%%%%%%%%%%%%%%%%%%%%%%%%%%%%%%%%%%%%%%%%%%%%%%%%%%%%
%%%%%%%%%%%%%%%%%%%%%%%%%%%%%%%%%%%%%%%%%%%%%%%%%%%%%%%%%%%%%%%%%%%%%%%%%%%%%%%%%%%%%%%%%%%%%
\bigskip
\begingroup
\makeatletter
\def\@biblabel#1{#1.}
\def\journal#1{\textit{\frenchspacing #1}}
\def\title#1{``#1''}
\def\booktitle#1{\textsl{#1}}
\def\v#1{\textbf{#1}}

\endgroup

%%% Local Variables:
%%% mode: latex
%%% TeX-master: "gibbs-ibvp"
%%% End:

\end{document}